\input amstex
\documentstyle{amsppt}
%
% Author:   Ruslan Abdulovich Sharipov
% Title:    Newtonian normal shift in multidimensional
%           Riemannian geometry
% Comment:  AmSTeX, 38 pages, amsppt style
% Address:  Rabochaya street 5, Ufa, 450003, Russia
% E-mail:   R_Sharipov@ic.bashedu.ru
%           ruslan-sharipov@usa.net
% Web-page: http://www.geocities.com/CapeCanaveral/Lab/5341
%
\nopagenumbers
\accentedsymbol\tx{\tilde x}
\accentedsymbol\ty{\tilde y}
\def\compos{\,\raise 1pt\hbox{$\sssize\circ$} \,}
\def\const{\operatorname{const}}
\def\negskp{\hskip -2pt}
\pagewidth{360pt}
\pageheight{606pt}
\rightheadtext{Newtonian normal shift \dots}
\topmatter
\title
NEWTONIAN NORMAL SHIFT\\IN MULTIDIMENSIONAL RIEMANNIAN GEOMETRY.
\endtitle
\author
Ruslan~A\.~Sharipov
\endauthor
\abstract
Explicit description for arbitrary Newtonian dynamical system
admitting the normal shift in Riemannian manifold of the dimension
$n\geqslant 3$ is found. On the base of this result the kinematics
of normal shift of hypersurfaces along trajectories of such system
is studied.
\endabstract
\address Rabochaya~str\.~5, 450003, Ufa, Russia
\endaddress
\email \vtop to 20pt{\hsize=280pt\noindent
R\_\hskip 1pt Sharipov\@ic.bashedu.ru\newline
ruslan-sharipov\@usa.net\vss}
\endemail
\urladdr
http:/\negskp/www.geocities.com/CapeCanaveral/Lab/5341
\endurladdr
\endtopmatter
\loadbold
\document
\head
1. Introduction.
\endhead
     In series of papers \cite{1--16} written in 1993--1996
a theory was constructed that determines and describes special
class of {\bf Newtonian dynamical systems admitting the normal
shift} of hypersurfaces in Riemannian and Finslerian manifolds.
On the base of these papers two theses were prepared: thesis
for the degree of {\bf Doctor of Sciences in Russia} \cite{17}
and thesis for the degree of {\bf Candidate of Sciences in Russia}
\cite{18}. However, some results included in thesis \cite{17} are
still not published in journals (see \cite{10} and \cite{11}).\par
     Moreover, when preparing thesis \cite{17}, in paper \cite{16}
an error was found. Eliminating this error led to new result that
consists in complete and exhausting description of all Newtonian
dynamical systems admitting the normal shift in
Riemannian\footnotemark\ manifolds of the dimension $n\geqslant 3$.
\footnotetext{This result has no direct generalization for the case
of dynamical systems in Finslerian ma\-nifolds (this case was considered
in \cite{10}, see also Chapter~\uppercase\expandafter{\romannumeral 8}
in thesis \cite{17}). But, nevertheless, it possibly has some analog
in Finslerian geometry. This problem is not yet studied.}The goal of
this paper is to explain this new result, and to give a description
for kinematics of normal shift of hypersurfaces, more detailed than
it was possible before now. On the base of the same result one can
get new (more simple) proof for the main theorem from unpublished
paper \cite{11} (see also \S\,7 in
Chapter~\uppercase\expandafter{\romannumeral 6} of thesis \cite{17}),
and one can answer the question by A.~V.~Bolsinov and A.~T.~Fomenko,
which they asked when author reported thesis \cite{17} in the seminar
of the Chair of Differential Geometry and its Applications at Moscow
State University.\adjustfootnotemark{-1}\par
     Classical construction of normal shift of hypersurfaces in
Riemannian manifold $M$ is well-known. In its original form it arises
in the case when $M$ is $\Bbb R^3$ with standard flat Euclidean
metric. Let $S$ be two-dimensional surface in $\Bbb R^3$. From each
point $p$ on $S$ we draw a segment of straight line in the direction
of normal vector $\bold n=\bold n(p)$. Denote by $p_t$ the second end
of this segment. When $p$ runs over $S$, point $p_t$ sweeps some other
surface $S_t$ as shown on Fig\. 1.1 below. So we have the map
$f_t\!:S\to S_t$ known as {\it classical normal shift} or as {\it Bonnet
transformation}.\par
\parshape 16 0pt 360pt 0pt 360pt 0pt 360pt 0pt 360pt 160pt 200pt
160pt 200pt 160pt 200pt 160pt 200pt 160pt 200pt 160pt 200pt
160pt 200pt 160pt 200pt 160pt 200pt 160pt 200pt 160pt 200pt 0pt 360pt
     When transferring from $M=\Bbb R^3$ to the case of arbitrary
Riemannian manifold $M$, we replace surfaces by hypersurfaces, and
straight line segments, connecting $p$ and $p_t$, by the segments
of geodesic lines. In this form classical construction of normal
shift of hypersurfaces is known as {\it geodesic normal shift}.
Construction of geodesic normal shift contains a numeric parameter
$t$, i\.~e\. it determines the whole family of hypersurfaces $S_t$.
\vadjust{\vskip -17pt\hbox to 0pt{\kern 5pt\hbox{\special{em:graph
Eds-01a.gif}}\hss}\vskip 17pt}When parameter $t$ is varied, point
$p_t$ moves along geodesic lines (here they are called {\it
trajectories of the shift}). In local coordinates $x^1,\,\ldots,\,x^n$
on $M$ they are described by ordinary differential equations
$$
\ddot x^k+\sum^n_{i=1}\sum^n_{j=1}\Gamma^k_{ij}\,\dot x^i\,\dot x^j=0,
\hskip -2em
\tag1.1
$$
where $k=1,\,\ldots,\,n$. By $\Gamma^k_{ij}=\Gamma^k_{ij}(x^1,
\ldots,x^n)$ in \thetag{1.1} we denote components of standard metric
connection $\Gamma$ for Riemannian metric $\bold g$ on $M$. The property
of normality of geodesic shift is expressed by the following well-known
fact.
\proclaim{Theorem 1.1} All hypersurfaces $S_t$ in the construction of
geodesic normal shift are perpendicular to trajectories of shift.
\endproclaim
In other words, trajectories of shift described by the equations
\thetag{1.1} cross each hypersurface $S_t$ transversally; at the
points of intersection they pass along normal vectors to $S_t$.
\par
    The idea of generalizing the construction of normal shift,
which was realized in papers \cite{1--16}, is very simple. It
consists in replacing \thetag{1.1} by slightly more complicated
ordinary differential equations in $t$:
$$
\ddot x^k+\sum^n_{i=1}\sum^n_{j=1}\Gamma^k_{ij}\,\dot x^i\,\dot x^j
=F^k(x^1,\ldots,x^n,\dot x^1,\ldots,\dot x^n).
\hskip -2em
\tag1.2
$$
When $M=\Bbb R^3$ and $\Gamma^k_{ij}=0$, the equations \thetag{1.2}
express Newton's second law: they describe the motion of a mass point
with unit mass in the force field determined by right hand sides of
these equations. In the case of arbitrary Riemannian manifold $M$
these equations, as appears, also have physical interpretation. They
describe the dynamics of complex mechanical systems with holonomic
constraints. Manifold $M$ arises as configuration space of such systems,
its dimension is determined by actual number of degrees of freedom
(upon resolving all constraints). Thereby $M$ is canonically equipped
with the structure of Riemannian manifold, its metric is given by
quadratic form of kinetic energy:
$$
K=\frac{1}{2}\sum^n_{i=1}\sum^n_{j=1}g_{ij}\,\dot x^i\,\dot x^j
$$
(see details in Chapter~\uppercase\expandafter{\romannumeral 2} of
thesis \cite{17}). Due to the analogy with Newton's second law the
equations \thetag{1.2} are called the equations {\it of Newtonian
dynamical system} on Riemannian manifold. Vector $\bold F$, whose
components are given by right hand sides of the equations \thetag{1.2},
is called a {\bf force vector}. It determines {\bf force field} of
Newtonian dynamical system \thetag{1.2}.\par
     Note that the choice of local coordinates in defining Newtonian
dynamical system \thetag{1.2} is of no importance. By the change
of local coordinates the shape of the equations remains unchanged,
thought the components of connections and components of force vector
are transformed according to standard formulas, which are well-known
from course of differential geometry (see \cite{19--22}). This property
expresses {\bf coordinate covariance} of differential equations
\thetag{1.2}.\par
\head
2. Newtonian normal shift of hypersurfaces.
\endhead
\parshape 24 0pt 360pt 0pt 360pt 0pt 360pt 0pt 360pt
0pt 360pt 0pt 360pt 0pt 360pt 0pt 360pt 0pt 360pt 0pt 360pt
0pt 360pt 160pt 200pt 160pt 200pt 160pt 200pt 160pt 200pt
160pt 200pt 160pt 200pt 160pt 200pt 160pt 200pt 160pt 200pt
160pt 200pt 160pt 200pt 160pt 200pt
0pt 360pt
     Having formulated the idea of generalizing the construction of
geodesic normal shift, we shall describe how it was realized in papers
\cite{1--16}. Let $S$ be some hypersurface in $M$ and let $p$ be some
point on $S$. In local coordinates $x^1,\,\ldots,\,x^n$ on $M$ such
point $p$ is characterized by its coordinates $x^1(p),\,\ldots,\,x^n(p)$
and by normal vector $\bold n(p)$ at this point. Let's use Newtonian
dynamical system \thetag{1.2} in order to define a shift of hypersurface
$S$. With this aim let's associate each point $p\in S$ with the 
following initial data for the system of differential equations
\thetag{1.2}:
$$
\xalignat 2
&\quad x^k\,\hbox{\vrule height 8pt depth 8pt width 0.5pt}_{\,t=0}
=x^k(p),
&&\dot x^k\,\hbox{\vrule height 8pt depth 8pt width 0.5pt}_{\,t=0}=
\nu(p)\cdot n^k(p).\hskip -2em
\tag2.1
\endxalignat
$$
Here $n^k(p)$ is $k$-th component of normal vector $\bold n(p)$, and
\vadjust{\vskip 16pt\hbox to 0pt{\kern 5pt\hbox{\special{em:graph
Eds-01b.gif}}\hss}
\vskip -16pt}$\nu(p)$ is some scalar quantity depending on the point $p\in S$.
Solving Cauchy problem with initial data \thetag{2.1} for the equations
\thetag{1.2}, we obtain a set of functions
$$
\cases x^1=x^1(t,p),\\ .\ .\ .\ .\ .\ .\ .\ .\ .\ \\
x^n=x^n(t,p).\endcases\hskip -1em
\tag2.2
$$
These functions define in parametric form the trajectory $r=r(t,p)$ of
Newtonian dynamical system with force field $\bold F$. This trajectory
at initial instant of time $t=0$ crosses hypersurface $S$ at the point
$p$, passing in the direction of unitary normal vector $\bold n(p)$.
Parameter $\nu(p)$ in initial data \thetag{2.1} determines the modulus
of initial velocity for this trajectory:
$$
\bold v\,\hbox{\vrule height 8pt depth 8pt width 0.5pt}_{\,t=0}=
\nu(p)\cdot\bold n(p).\hskip -2em
\tag2.3
$$
Choice of local coordinates in defining trajectory $r=r(t,p)$ is of
no importance. Change of local coordinates changes functions \thetag{2.2},
but it doesn't change the curve $r=r(t,p)$. This is due to coordinate
covariance of differential equations \thetag{1.2} and coordinate
covariance of initial data \thetag{2.1}.\par
    Drawing trajectories $r=r(t,p)$ outgoing from all points $p\in S$
and taking points $p_t=r(t,p)$ that corresponds to some fixed value
of parameter $t$, we obtain the hypersurface $S_t$ and displacement map
$f_t\!:S\to S_t$. However, we should remember two nuances. Parameter $t$
for trajectory $r=r(t,p)$ of dynamical system \thetag{1.2} do not
coincide with its length. The range of this parameter $t$ always includes
initial point $t=0$, but it can be a restricted interval
$$
t_1(p)<t<t_2(p).
$$
Upper and lower bounds of this interval in general case depend on the
point $p\in S$. Hence for a fixed value of $t$ the displacement map
$f_t\!:S\to S_t$ can be defined not for all points $p\in S$.\par
    Second nuance is due to singular points (caustics) that may appear
on the hypersurface $S_t$ for large enough values of parameter $t$.
This imposes one more restriction onto the range of parameter $t$.
Note that this restriction is present in classical construction of
geodesic normal shift as well.\par
    The above two nuances restrict possible range of parameter $t$.
However, if we are interested in small values of $t$ only (as below),
we can use the following lemma.
\proclaim{Lemma 2.1} If parameter $\nu(p)$ in \thetag{2.3} is a smooth
nonzero function on the hypersurface $S$, then for each $p\in S$ there
exists some neighborhood $S'=O_{\!S}(p)$ on $S$ and there exists
a number $\varepsilon$ such that displacement maps $f_t\!:S'\to S'_t$
are defined for all $t\in(-\varepsilon,\,+\varepsilon)$. They form
smooth one-parametric family of diffeomorphisms.
\endproclaim
   This lemma is an immediate consequence of theorem on existence,
uniqueness, and smooth dependence of initial data for the solution
of Cauchy problem for systems of ODE's (see \cite{23} and \cite{24}).
Taking into account lemma~2.1, we can consider displacement maps
$f_t\!:S\to S_t$, which are possibly defined only locally on $S$, as
a construction of {\bf shift} of hypersurface $S$ {\bf along trajectories}
of Newtonian dynamical system with force field $\bold F$. Function $\nu(p)$
on $S$ is a parameter in such construction of shift.\par
     Shift of hypersurface $S$ along trajectories of dynamical system
\thetag{1.2}, as it was constructed above, possess the property of
normality at the initial instant of time $t=0$. This means that trajectories
of shift are passing through initial hypersurface $S$ along normal vectors
on it. Does this property persist for $t\neq 0$, i\.~e\. can we prove
theorem similar to theorem~1.1\,? The answer to this question in general
case is negative (see examples in \cite{18}). But there are special cases,
when the property of normality persists for all instants of time. We
describe them by formulating the following definition.
\definition{Definition 2.1} Shift $f_t\!:S\to S_t$ of hypersurface $S$
along trajectories of Newtonian dynamical system with force field
$\bold F$ is called a {\bf normal shift} if all hypersurfaces $S_t$
(for all permissible values of parameter $t$) are orthogonal to
trajectories of shift.
\enddefinition
\head
3. Dynamical systems admitting the normal shift.
\endhead
     What does the property of normality for shift $f_t\!:S\to S_t$ depend
on\,? On the choice of hypersurface $S$\,? On the force field $\bold F$ of
Newtonian dynamical system\,? We also have the opportunity to choose the
function $\nu(p)$ on $S$. In the case of identically zero force field
$\bold F=0$ (which corresponds to geodesic flow on $M$) the choice
$\nu(p)=1$ provides normality condition from definition~2.1 for arbitrary
initial hypersurface $S\subset M$. Are there some other force fields with
similar property\,? The aim to know this was the motivation for writing
preprint \cite{1}. In this preprint we introduced the concept of {\bf
Newtonian dynamical system admitting the normal shift}. This concept has
become a central point for later investigations.
\definition{Definition 3.1} Newtonian dynamical system on Riemannian
manifold $M$ is called a system {\bf admitting the normal shift} if
for any hypersurface $S$ in $M$, and for any point $p_0\in S$, there is
a neighborhood $S'=O_{\!S}(p_0)$ of the point $p_0$ on $S$, and there is
a smooth function $\nu(p)$ in $S'$, such that the shift $f_t\!:S'\to S'_t$
defined by the function $\nu(p)$ is a normal shift along trajectories of
considered dynamical in the sense of definition~2.1.
\enddefinition
    The condition stated in definition~3.1 was called {\bf the normality
condition} for Newtonian dynamical system with force field $\bold F$.
First we considered the case $M=\Bbb R^2$ (see \cite{1}). In \cite{1}
we derived partial differential equations for the components of force
field $\bold F$ which, when being fulfilled, are sufficient to provide
normality condition from definition~3.1. These equations were called
the {\bf normality equations} or, more exactly, {\bf weak normality
equations}. In preprint \cite{1} we also constructed first non-trivial
examples of dynamical systems that admit normal shift. When generalizing
these results from $M=\Bbb R^2$ to multidimensional case $M=\Bbb R^n$ in
preprint \cite{1}, we have found that weak normality equations should
be supplemented by so called additional normality equations. All above
results in brief form were announced in \cite{4}. Their full version
were published in \cite{2} and \cite{3}. Later in papers \cite{6} and
\cite{7} they were generalized for the case of arbitrary Riemannian
manifold $M$. Main result of paper \cite{6} is the derivation of weak
normality equations for this more complicated geometric situation. We
write these equations without comments so far:
$$
\cases
\dsize\sum^n_{i=1}\left(v^{-1}\,F_i+\shave{\sum^n_{j=1}}
\tilde\nabla_i\left(N^j\,F_j\right)\right)P^i_k=0,\\
\vspace{2.5ex}
\aligned
 &\sum^n_{i=1}\sum^n_{j=1}\left(\nabla_iF_j+\nabla_jF_i-2\,v^{-1}
 \,F_i\,F_j\right)N^j\,P^i_k\,+\\
 \vspace{0.5ex}
 &+\sum^n_{i=1}\sum^n_{j=1}\left(\frac{F^j\,\tilde\nabla_jF_i}{v}
 -\sum^n_{r=1}\frac{N^r\,N^j\,\tilde\nabla_jF_r}{v}\,F_i\right)P^i_k=0.
 \endaligned
\endcases\hskip -3em
\tag3.1
$$
Additional normality equations for the force field $\bold F$ of
Newtonian dynamical systems on Riemannian manifolds were derived
in \cite{7}. They look like 
$$
\cases
\aligned
&\sum^n_{i=1}\sum^n_{j=1}P^i_\varepsilon\,P^j_\sigma\left(
\,\shave{\sum^n_{m=1}}N^m\,\frac{F_i\,\tilde\nabla_mF_j}{v}
-\nabla_iF_j\right)=\\
\vspace{0.5ex}
&\quad=\sum^n_{i=1}\sum^n_{j=1}P^i_\varepsilon\,P^j_\sigma\left(
\,\shave{\sum^n_{m=1}}N^m\,\frac{F_j\,\tilde\nabla_mF_i}{v}
-\nabla_jF_i\right),
\endaligned\\
\vspace{2.5ex}
\dsize\sum^n_{i=1}\sum^n_{j=1}P^j_\sigma\,\tilde\nabla_jF^i\,
P^\varepsilon_i=\sum^n_{i=1}\sum^n_{j=1}\sum^n_{m=1}
\frac{P^j_m\,\tilde\nabla_jF^i\,P^m_i}{n-1}\,P^\varepsilon_\sigma.
\endcases\hskip -3em
\tag3.2
$$
Results of papers \cite{6} and \cite{7} were announced in \cite{12}.
In deriving the equations \thetag{3.2} we have found that it is
convenient to make slight modification of definition~3.2. We added
normalizing condition for the function $\nu(p)$:
$$
\nu(p_0)=\nu_0.\hskip -2em
\tag3.3
$$
As a result we obtained {\bf strong normality condition}. It is
formulated as follows.
\definition{Definition 3.2} Newtonian dynamical system on Riemannian
manifold $M$ is called a system admitting the normal shift {\bf in
strong sense} if for any hypersurface $S$ in $M$, for any point
$p_0\in S$, and for any real number $\nu_0\neq 0$ there is a neighborhood
$S'=O_{\!S}(p_0)$ of the point $p_0$ on $S$, and there is a smooth
function $\nu(p)$ in $S'$ normalized by the condition \thetag{3.3}, such
that the shift $f_t\!:S'\to S'_t$ defined by the function $\nu(p)$ is
a normal shift along trajectories of considered dynamical in the sense
of the definition~2.1.
\enddefinition
    Strong normality condition from definition~3.2 implies normality
condition formulated in definition~3.1. The relation of strong normality
condition with normality equations \thetag{3.1} and \thetag{3.2} is
described by the following theorem.
\proclaim{Theorem 3.1} Newtonian dynamical system on Riemannian manifold
$M$ admits the normal shift in strong sense if and only if its force
field $\bold F$ satisfies normality equations \thetag{3.1} and \thetag{3.2}
simultaneously.
\endproclaim
     Theorem~3.1 was proved in \cite{13}. Detailed version of this proof
can be found in thesis \cite{17}. We shall not give this proof here,
since this would require to reproduce many details of derivation of the
normality equations \thetag{3.1} and \thetag{3.2}, and hence this would
be doubling for the papers \cite{6} and \cite{7}. Instead, we shall give
comments to normality equations \thetag{3.1} and \thetag{3.2}.
\head
4. Extended algebra of tensor fields.
\endhead
     Normality equations \thetag{3.1} and \thetag{3.2}, as well as
the  equations  \thetag{1.2},  possess  the  property   of   {\bf 
coordinate
covariance}. Writing these equations, we assume that some local
coordinates $x^1,\,\ldots,\,x^n$ in $M$ are chosen. Under the change
of local coordinates all quantities, which are contained in the equations
\thetag{3.1} and \thetag{3.2}, do change according some definite rules.
However, this do not change the shape of these equations. Components
of force vector $\bold F$ in \thetag{3.1} and \thetag{3.2} depend on
$x^1,\,\ldots,\,x^n$, and on components of velocity vector $\bold v$;
the latter is a tangent vector at the point $p$ with coordinates $x^1,\,
\ldots,\,x^n$. This means that vector $\bold F$ depend on the point
$q=(p,\bold v)$ of the tangent bundle $TM$.
\definition{Definition 4.1} Vector-function $\bold F$ that for each
point $q=(p,\bold v)$ of tangent bundle $TM$ puts into the correspondence
some vector from tangent space $T_p(M)$ at the point $p=\pi(q)$ is called
an {\bf extended vector field} on the manifold $M$.
\enddefinition
Here $\pi\!:TM\to M$ is a map of canonical projection from $TM$ to
the base manifold $M$. Let's consider the following tensor product:
$$
T^r_s(p,M)=\overbrace{T_p(M)\otimes\ldots\otimes T_p(M)}^{\text{$r$
times}}\otimes\underbrace{T^*_p(M)\otimes\ldots\otimes T^*_p(M)}_{\text{$s$
times}}.
$$
Linear space $T^r_s(p,M)$ is called a space of tensors of the type $(r,s)$
at the point $p$ of the manifold $M$. Elements of this space are called
{\it $r$-times contravariant and $s$-times covariant tensors}, or simply
tensors of type $(r,s)$ at the point $p\in M$.
\definition{Definition 4.2} Tensor-valued function $\bold X$ that for
each point $q$ of tangent bundle $TM$ puts into the correspondence some
tensor from the space $T^r_s(p,M)$ at the point $p=\pi(q)$ is called an
{\bf extended tensor field} of type $(r,s)$ on the manifold $M$.
\enddefinition
    In local coordinates $x^1,\,\ldots,\,x^n$ on the manifold $M$ extended
tensor fields are expressed by the functions of double-set of arguments:
$$
X^{i_1\ldots\,i_r}_{j_1\ldots\,j_s}(x^1,\ldots,x^n,v^1,\ldots,v^n).
\hskip -2em
\tag4.1
$$
In the normality equations \thetag{3.1} and \thetag{3.2} we can see
components of several extended tensor fields. Velocity vector $\bold v$
by itself can be considered as extended vector field on $M$. Its
modulus $v=|\bold v|$ is an extended scalar. Extended vector field
$\bold N$ with components $N^1,\,\ldots,\,N^n$, which are contained
in the equations \thetag{3.1} and \thetag{3.2}, is defined as the
following quotient:
$$
\bold N=\frac{\bold v}{v}=\frac{\bold v}{|\bold v|}.\hskip -2em
\tag4.2
$$
This is the field of unitary vectors directed along the vector of
velocity. And finally, in equations \thetag{3.1} and \thetag{3.2}
we have components of operator field $\bold P$. This is the field
of orthogonal projectors onto the hyperplane perpendicular to the
velocity vector. Components $P^i_j$ of this field are given by
the formula
$$
P^i_j=\delta^i_j-N^i\,N_j.\hskip -2em
\tag4.3
$$\par
    Let $T^r_s(M)$ be the set of smooth extended tensor fields of
type $(r,s)$. This set has the structure of module over the ring
of extended scalar fields. The following sum is a graded algebra
over this ring with respect to tensorial multiplication:
$$
\bold T(M)=\bigoplus^\infty_{r=0}\bigoplus^\infty_{s=0}T^r_s(M).
\hskip -2em
\tag4.4
$$
Algebra \thetag{4.4} is called an {\bf extended algebra of tensor
fields} on the manifold $M$. In extended algebra of tensor fields
$\bold T(M)$ one can define two operations of covariant differentiation,
we denote them by $\nabla$ and $\tilde\nabla$:
$$
\xalignat 2
&\nabla\!: T^r_s(M)\to T^r_{s+1}(M),
&&\tilde\nabla\!: T^r_s(M)\to T^r_{s+1}(M).
\endxalignat
$$
In local coordinates the result of applying covariant differentiation
$\nabla$ to a tensor field $\bold X$ with components \thetag{4.1} is
expressed by the following formula:
$$
\aligned
\nabla_m&X^{i_1\ldots\,i_r}_{j_1\ldots\,j_s}=\frac{\partial
X^{i_1\ldots\,i_r}_{j_1\ldots\,j_s}}{\partial x^m}
-\sum^n_{a=1}\sum^n_{b=1}v^a\,\Gamma^b_{ma}\,\frac{\partial
X^{i_1\ldots\,i_r}_{j_1\ldots\,j_s}}{\partial v^b}\,+\\
&+\sum^r_{k=1}\sum^n_{a_k=1}\!\Gamma^{i_k}_{m\,a_k}\,X^{i_1\ldots\,
a_k\ldots\,i_r}_{j_1\ldots\,\ldots\,\ldots\,j_s}
-\sum^s_{k=1}\sum^n_{b_k=1}\!\Gamma^{b_k}_{m\,j_k}
X^{i_1\ldots\,\ldots\,\ldots\,i_r}_{j_1\ldots\,b_k\ldots\,j_s}.
\endaligned\hskip -2em
\tag4.5
$$
The result of applying $\tilde\nabla$ to $\bold X$ is expressed by less
complicated formula:
$$
\tilde\nabla_mX^{i_1\ldots\,i_r}_{j_1\ldots\,j_s}=\frac{\partial
X^{i_1\ldots\,i_r}_{j_1\ldots\,j_s}}{\partial v^m}.
\hskip -2em
\tag4.6
$$
Formula \thetag{4.6} for the components of the field $\tilde\nabla\bold X$
contains \pagebreak only derivatives with respect to components of velocity
vector. Therefore $\tilde\nabla$ is called {\bf velocity gradient}.
Covariant differentiation $\nabla$ defined by formula \thetag{4.5} is called
{\bf spatial gradient}.\par
    Defining operators $\nabla$ and $\tilde\nabla$ by means of formulas
\thetag{4.5} and \thetag{4.6}, we assume that some local coordinates are
chosen. This way of defining $\nabla$ and $\tilde\nabla$ is quite
sufficient for our purposes in the theory of newtonian dynamical systems
admitting the normal shift. But there is another (invariant) way of
defining these operators. It is based on the analysis of differentiations
in the extended algebra of tensor fields.
\definition{Definition 4.3} The map $D\!:\bold T(M)\to\bold T(M)$ is
called a {\bf differentiation} in extended algebra of tensor fields
if the following conditions are fulfilled:
\roster
\rosteritemwd=10pt
\item compatibility with grading: $D(T^r_s(M))\subset
      T^r_s(M)$;
\item $\Bbb R$-linearity: $D(\bold X+\bold Y)=D(\bold X)
      +D(\bold Y)$ and $D(\lambda\bold X)=\lambda D(\bold X)$
      for $\lambda\in\Bbb R$;
\item permutability with contractions: $D(C(\bold X))=
      C(D(\bold X))$;
\item Leibniz rule: $D(\bold X\otimes\bold Y)=D(\bold X)
      \otimes\bold Y+\bold X\otimes D(\bold Y)$.
\endroster
\enddefinition
    Among the results of thesis \cite{17} the following structural theorem
is worth to mention here. It describes the structure of all differentiations
in extended algebra of tensor fields $\bold T(M)$.
\proclaim{Theorem 4.1} Let $M$ be smooth real manifold equipped with some
extended affine connection $\Gamma$. Then each differentiation $D$ in
extended algebra of tensor fields $\bold T(M)$ on this manifold breaks into
the sum
$$
D=\nabla_{\bold X}+\tilde\nabla_{\bold Y}+\bold S,
$$
where $\bold X$ and $\bold Y$ are some extended vector fields, and $\bold S$
is a degenerate differentiation defined by some extended tensor field $\bold
S$ of type $(1,1)$ in $M$.
\endproclaim
    Theorem~4.1 is an analog of structural theorem for differentiations in
the algebra of ordinary (not extended) tensor fields (see \cite{19}).
\head
5. Reduction of normality equations in the dimension $n\geqslant 3$.
\endhead
    If we take into account formulas \thetag{4.2}, \thetag{4.3},
\thetag{4.5}, and \thetag{4.6}, we see that normality equations
\thetag{3.1} and \thetag{3.2} form strongly overdetermined system
of partial differential equations with respect to components of
force field $\bold F$ of Newtonian dynamical system. Analysis of this
system (see \cite{16}), is based on scalar ansatz
$$
F_k=A\,N_k-|\bold v|\,\sum^n_{i=1}P^i_k\,\tilde\nabla_i A.
\hskip -2em
\tag5.1
$$
Formula \thetag{5.1} expresses components of force vector through one
scalar field $A$, which is interpreted as the projection of $\bold F$
onto the direction of velocity vector. This formula follows from first
part of equations in the system \thetag{3.1}. Therefore, when substituting
\thetag{5.1} into weak normality equations \thetag{3.1}, first part of
these equations appears to be identically fulfilled. While second part
is brought to
$$
\aligned
\sum^n_{s=1}&\left(\nabla_sA+|\bold v|\sum^n_{q=1}\sum^n_{r=1}
P^{qr}\,\tilde\nabla_qA\,\tilde\nabla_r\tilde\nabla_sA\right.-\\
&-\left.\shave{\sum^n_{r=1}}N^r\,A\,\tilde\nabla_r\tilde\nabla_sA
-|\bold v|\shave{\sum^n_{r=1}}N^r\,\nabla_r\tilde\nabla_sA\right)
P^s_k=0.
\endaligned\hskip -2em
\tag5.2
$$
For $n=2$ the equations \thetag{5.2} exhaust whole list of reduced
normality equations. The matter is that in two-dimensional case, as
we mentioned above, additional normality equations do not arise at
all. While the equations \thetag{5.2} are reduced to the only one
nonlinear partial differential equation for the function $A(x^1,x^2,
v^1,v^2)$. Detailed study of this equation is given in the thesis
by A\.~Yu\.~Boldin \cite{18}.\par
     In multidimensional case $n\geqslant 3$ the process of reducing
normality equations can be moved much further. Substituting \thetag{5.1}
into the first part of additional normality equations brings them to
the following form:
$$
\aligned
\sum^n_{s=1}&\sum^n_{r=1}P^r_\sigma\,P^s_\varepsilon\left(
\nabla_r\tilde\nabla_sA+\shave{\sum^n_{q=1}}\tilde\nabla_rA\,N^q\,
\tilde\nabla_q\tilde\nabla_sA\right)=\\
&=\sum^n_{s=1}\sum^n_{r=1}P^r_\sigma\,P^s_\varepsilon\left(
\nabla_s\tilde\nabla_rA+\shave{\sum^n_{q=1}}\tilde\nabla_sA\,N^q\,
\tilde\nabla_q\tilde\nabla_rA\right).
\endaligned\hskip -2em
\tag5.3
$$
Similarly, substituting \thetag{5.1} into second part of the equations
\thetag{3.2} gives
$$
\sum^n_{r=1}\sum^n_{s=1}P^r_\sigma\,\tilde\nabla_r\tilde\nabla_sA\,
P^{s\varepsilon}=\lambda\,P^\varepsilon_\sigma.\hskip -2em
\tag5.4
$$
Here $\lambda$ is a scalar quantity, the value of which is uniquely
determined by the equations \thetag{5.4} even if we do not know it
a priori:
$$
\lambda=\sum^n_{r=1}\sum^n_{s=1}\frac{P^{rs}\,\tilde\nabla_r
\tilde\nabla_sA}{n-1}.
$$
The equations \thetag{5.4} are most remarkable. According to the
formula \thetag{4.6} they contain only the derivatives with respect
to the variables $v^1,\,\ldots,\,v^n$. This corresponds to varying
the function $A(x^1,\ldots,x^n,v^1,\ldots,v^n)$ within the fiber of
$TM$ over the fixed point with coordinates $x^1,\,\ldots,\,x^n$ in
the base manifold $M$.
\definition{Definition 5.1} Extended tensor field $\bold X$ on Riemannian
manifold $M$ is called {\bf fiberwise spherically symmetric} if $|\bold v_1|
=|\bold v_2|$ implies $\bold X(p,\bold v_1)=\bold X(p,\bold v_2)$.
\enddefinition
In other words, fiberwise spherically symmetric extended tensor fields
depend only on modulus of velocity vector within fibers of tangent
bundle $TM$. Such fields naturally arises in the analysis of the equations
\thetag{5.4}. Here we have the following theorem proved in paper
\cite{16}.
\proclaim{Theorem 5.1} Extended Scalar field $A$ on Riemannian manifold
$M$ satisfies equations \thetag{5.4} if and only if it is given by
formula
$$
A=a+\sum^n_{i=1}b_i\,v^i,\hskip -2em
\tag5.5
$$
where $a$ is some fiberwise spherically symmetric scalar field, and
$b_i$ are components of some fiberwise spherically symmetric covectorial
field $\bold b$.
\endproclaim
     Further substitution of \thetag{5.5} into the equations \thetag{5.2}
and \thetag{5.3} yields the following equations with respect to fields
$a$ and $\bold b$:
$$
\gather
\left(\frac{\partial}{\partial x^s}+b_s\,\frac{\partial}{\partial v}
\right)a=\left(a\,\frac{\partial}{\partial v}\right)b_s,\hskip -2em
\tag5.6\\
\displaybreak
\left(\frac{\partial}{\partial x^s}+b_s\,\frac{\partial}{\partial v}
\right)b_r=\left(\frac{\partial}{\partial x^r}+b_r\,\frac{\partial}
{\partial v}\right)b_s.\hskip -2em
\tag5.7
\endgather
$$
Here $a=a(x^1,\ldots,x^n,v)$ and $b_i=b_i(x^1,\ldots,x^n,v)$. Variable
$v$ denotes the modulus of velocity vector: $v=|\bold v|$.\par
    Equations \thetag{5.6} and \thetag{5.7} is bound with most dramatic
instant in constructing theory of dynamical systems admitting the normal
shift. By deriving these equations in paper \cite{16} the mistake was made,
nonlinear terms in \thetag{5.7} were omitted. As a result \thetag{5.7}
looked like $\partial b_r/\partial x^s=\partial b_s/\partial x^r$. Further
analysis of erroneous equations has led to the ordinary differential equation
$$
y''=H_y(y'+1)+H_x\,\hskip -2em
\tag5.8
$$
where $H=H(x,y)$, $H_y=\partial H/\partial y$, $H_x=\partial H/\partial
x$. With the aim to find as more functions $H(x,y)$, for which the equation
\thetag{5.8} is explicitly solvable, as possible we considered the following
change of variables:
$$
\cases
\tx=\tx(x,y),\\ \ty=\ty(x,y).
\endcases\hskip -2em
\tag5.9
$$
The equations that could be brought to the form \thetag{5.8} by means of
change of variables \thetag{5.9} belong to the following class of equations:
$$
y''=P(x,y)+3\,Q(x,y)\,y'+3\,R(x,y)\,{y'}^2+S(x,y)\,{y'}^3.
\hskip -3em
\tag5.10
$$
Study of point transformations in the class of equations \thetag{5.10}
has the long history (see \cite{25--46}). However, we couldn't find an
answer to the question: how to extract the equations \thetag{5.10} that
could be brought to the form \thetag{5.8} by means of point transformation
\thetag{5.9}. This stimulated our own investigations (see \cite{47--50}).
We managed to get some results in describing classes of point equivalence
for the equations \thetag{5.10}. But now, since the error in \cite{16} is
found, these results have separate value, which is not related to the
theory of dynamical systems admitting the normal shift. And we are to
return to the equations \thetag{5.6} and \thetag{5.7}.\par
\head
6. Derivation of reduced normality equations.
\endhead
     For the beginning let's derive the \thetag{5.6} and \thetag{5.7}
by substituting \thetag{5.5} into the equations \thetag{5.2} and
\thetag{5.3}. Denote by $a'$ and $b'_i$ the derivatives
$$
\xalignat 2
&a'=\frac{\partial a}{\partial v},
&&b'_i=\frac{\partial b_i}{\partial v}.
\endxalignat
$$
Let's do the calculations necessary for substituting \thetag{5.5} into
\thetag{5.2} and \thetag{5.3}:
$$
\allowdisplaybreaks
\gather
\nabla_sA=\nabla_sa+\sum^n_{i=1}\nabla_sb_i\,v^i,
\hskip -2em
\tag6.1\\
\nabla_r\tilde\nabla_sA=\left(\nabla_ra'+\shave{\sum^n_{i=1}}
\nabla_rb'_i\,v^i\right)N_s+\nabla_rb_s.\hskip -2em
\tag6.2\\
\tilde\nabla_sA=\left(a'+\shave{\sum^n_{i=1}}b'_i\,v^i\right)
N_s+b_s,\hskip -2em
\tag6.3\\
\aligned
\tilde\nabla_r&\tilde\nabla_sA=\left(a''+\shave{\sum^n_{i=1}}b''_i
\,v^i\right)N_r\,N_s\,+\\
&+\,b'_s\,N_r+b'_r\,N_s+\left(\frac{a'}{v}+\shave{\sum^n_{i=1}}b'_i
\,N^i\right)P_{rs}.
\endaligned\hskip -3em
\tag6.4
\endgather
$$
From formulas \thetag{6.3} and \thetag{6.4} for derivatives we obtain
the following relations:
$$
\gather
\sum^n_{r=1}P^r_\sigma\,\tilde\nabla_rA=\sum^n_{r=1}P^r_\sigma\,b_r,
\hskip -2em
\tag6.5\\
\sum^n_{s=1}\sum^n_{q=1}P^s_\varepsilon\,N^q\,\tilde\nabla_q
\,\tilde\nabla_sA=\sum^n_{s=1}P^s_\varepsilon\,b'_s.\hskip -2em
\tag6.6
\endgather
$$
Let's combine \thetag{6.5} and \thetag{6.6}. As a result we get the
relationship 
$$
\sum^n_{s=1}\sum^n_{r=1}\sum^n_{q=1}P^r_\sigma\,P^s_\varepsilon
\,\tilde\nabla_rA\,N^q\,\tilde\nabla_q\tilde\nabla_sA=\sum^n_{s=1}
\sum^n_{r=1}P^r_\sigma\,P^s_\varepsilon\,b_r\,b'_s.
$$
Then let's multiply both sides of the relationship \thetag{6.2} by
$P^r_\sigma\,P^s_\varepsilon$ and contract with respect to pair of
indices $r$ and $s$. This yields one more relationship:
$$
\sum^n_{s=1}\sum^n_{r=1}P^r_\sigma\,P^s_\varepsilon\,
\nabla_r\tilde\nabla_sA=\sum^n_{s=1}\sum^n_{r=1}P^r_\sigma\,
P^s_\varepsilon\,\nabla_rb_s.
$$
Now, if we add above two relationships, we get the result of substituting
\thetag{5.5} into the left hand side of the equation \thetag{5.3}:
$$
\aligned
\sum^n_{s=1}&\sum^n_{r=1}P^r_\sigma\,P^s_\varepsilon\left(
\nabla_r\tilde\nabla_sA+\shave{\sum^n_{q=1}}\tilde\nabla_rA\,
N^q\,\tilde\nabla_q\tilde\nabla_sA\right)=\\
&=\sum^n_{s=1}\sum^n_{r=1}P^r_\sigma\,P^s_\varepsilon\,\nabla_rb_s
+\sum^n_{s=1}\sum^n_{r=1}P^r_\sigma\,P^s_\varepsilon\,b_r\,b'_s.
\endaligned\hskip -2em
\tag6.7
$$
Similarly we calculate right hand side of the equation \thetag{5.3}:
$$
\aligned
\sum^n_{s=1}&\sum^n_{r=1}P^r_\sigma\,P^s_\varepsilon\left(
\nabla_s\tilde\nabla_rA+\shave{\sum^n_{q=1}}\tilde\nabla_sA\,
N^q\,\tilde\nabla_q\tilde\nabla_rA\right)=\\
&=\sum^n_{s=1}\sum^n_{r=1}P^r_\sigma\,P^s_\varepsilon\,\nabla_sb_r
+\sum^n_{s=1}\sum^n_{r=1}P^r_\sigma\,P^s_\varepsilon\,b_s\,b'_r.
\endaligned\hskip -2em
\tag6.8
$$
On the base of \thetag{6.7} and \thetag{6.8} we conclude that the
equation \thetag{5.3} is reduced to
$$
\pagebreak
\sum^n_{s=1}\sum^n_{r=1}P^r_\sigma\,P^s_\varepsilon\,\bigl(\nabla_rb_s
+b_r\,b'_s-\nabla_sb_r-b_s\,b'_r\bigr)=0.\hskip -2em
\tag6.9
$$
For the further analysis of the obtained equations \thetag{6.9} one
should use the peculiarity of covectorial field $\bold b$ from extended
algebra of tensor fields on $M$. Components of this field $b_1,\,\ldots,
\,b_n$ depend only on modulus of velocity vector, but they do not depend
on its direction. By calculating derivatives $\nabla_rb_s$ and
$\nabla_sb_r$ in \thetag{6.9} we apply the following theorem.
\proclaim{Theorem 6.1} Let $X^{i_1\ldots\,i_r}_{j_1\ldots\,j_s}=
X^{i_1\ldots\,i_r}_{j_1\ldots\,j_s}(x^1,\ldots,x^n,v)$ be components
of fiberwise spherically symmetric tensor field $\bold X$ from extended
algebra $\bold T(M)$. Then components of spatial gradient $\nabla\bold X$
for this field are given by formula
$$
\aligned
\nabla_mX^{i_1\ldots\,i_r}_{j_1\ldots\,j_s}=&\frac{\partial
X^{i_1\ldots\,i_r}_{j_1\ldots\,j_s}}{\partial x^m}+
\sum^r_{k=1}\sum^n_{a_k=1}\!\Gamma^{i_k}_{m\,a_k}\,X^{i_1\ldots\,
a_k\ldots\,i_r}_{j_1\ldots\,\ldots\,\ldots\,j_s}\,-\\
&-\sum^s_{k=1}\sum^n_{b_k=1}\!\Gamma^{b_k}_{m\,j_k}
X^{i_1\ldots\,\ldots\,\ldots\,i_r}_{j_1\ldots\,b_k\ldots\,j_s}.
\endaligned\hskip -3em
\tag6.10
$$
\endproclaim
\demo{Proof} Formula \thetag{6.10} is obtained as a result of reduction
from formula \thetag{4.5}. For the components of fiberwise spherically
symmetric tensor field $\bold X$ (as in the statement of theorem) natural
arguments are $x^1,\,\ldots,\,x^n$, and $v$, where
$$
v=|\bold v|=\sqrt{\shave{\sum^n_{i=1}\sum^n_{j=1}}g_{ij}(x^1,\ldots,
x^n)\,v^i\,v^j\,}.\hskip -3em
\tag6.11
$$
While partial derivatives in formula \thetag{4.5} are assumed to be
respective to the variables $x^1,\,\ldots,\,x^n,\,v^1,\,\ldots,\,v^n$.
Recalculation of these derivatives to natural variables for spherically
symmetric field consists in the following substitutions:
$$
\gather
\frac{\partial X^{i_1\ldots\,i_r}_{j_1\ldots\,j_s}}{\partial v^b}
\text{ \ by \ }\frac{\partial X^{i_1\ldots\,i_r}_{j_1\ldots\,j_s}}
{\partial v}\cdot\frac{\partial v}{\partial v^b},\hskip -3em
\tag6.12\\
\frac{\partial X^{i_1\ldots\,i_r}_{j_1\ldots\,j_s}}{\partial x^m}
\text{ \ by \ }\frac{\partial X^{i_1\ldots\,i_r}_{j_1\ldots\,j_s}}
{\partial x^m}+\frac{\partial X^{i_1\ldots\,i_r}_{j_1\ldots\,j_s}}
{\partial v}\cdot\frac{\partial\,v}{\partial x^m}.\hskip -3em
\tag6.13
\endgather
$$
Derivatives $\partial\,v/\partial v^b$ and $\partial\,v/\partial x^m$
are calculated due to \thetag{6.11}. Upon finding explicit expressions
for these derivatives and upon making substitutions \thetag{6.12} and
\thetag{6.13} in formula \thetag{4.5}, we get two extra summands:
$$
\sum^n_{a=1}\sum^n_{b=1}\frac{X^{i_1\ldots\,i_r}_{j_1\ldots\,j_s}}
{\partial v}\,\frac{1}{2}\,\frac{\partial g_{ab}}{\partial x^m}\,
\frac{v^a\,v^b}{v}\text{ \ \ and \ \ }
-\sum^n_{a=1}\sum^n_{b=1}v^a\,\Gamma^b_{ma}\,\frac{\partial
X^{i_1\ldots\,i_r}_{j_1\ldots\,j_s}}{\partial v}\,N_b.
$$
If we take into account the explicit formula for components of metric
connection
$$
\Gamma^b_{ma}=\frac{1}{2}\sum^n_{q=1}g^{bq}\left(\frac{\partial g_{qa}}
{\partial x^m}+\frac{\partial g_{mq}}{\partial x^a}-\frac{\partial g_{ma}}
{\partial x^q}\right)\hskip -3em
\tag6.14
$$
(see \cite{19--22}), then we easily see that above two summands cancel
each other. As a result formula \thetag{4.5} transforms into the form
\thetag{6.10}.
\qed\enddemo
\proclaim{Corollary} Components $\omega_{rs}=\nabla_rb_s+b_r\,b'_s
-\nabla_sb_r-b_s\,b'_r$ of skew-symmetric extended tensor field
$\boldsymbol\omega$ in the equations \thetag{6.9} depend only on
modulus of velocity vector $|\bold v|$, but they do not depend on
the direction of the vector $\bold v$.
\endproclaim
     This fact immediately follows from the formula \thetag{6.10}.
It allows us to make further simplifications in the equations
\thetag{6.9}. Let $\bold c$ and $\bold d$ be arbitrary two vectors
from tangent space $T_p(M)$. In multidimensional case $n\geqslant 3$
we can rotate velocity vector $\bold v$, keeping its modulus unchanged,
and can direct it so that it will be perpendicular to vectors $\bold c$
and $\bold d$ simultaneously. Then
$$
\xalignat 2
&\sum^n_{\sigma=1}P^r_\sigma\,c^{\sigma}=c^r,
&&\sum^n_{\varepsilon=1}P^s_\varepsilon\,d^{\varepsilon}=d^s.
\endxalignat
$$
Therefore the equations \thetag{6.9} are transformed as follows:
$$
\sum^n_{s=1}\sum^n_{r=1}c^r\,d^s\,\bigl(\nabla_rb_s
+b_r\,b'_s-\nabla_sb_r-b_s\,b'_r\bigr)=0.\hskip -3em
\tag6.15
$$
Since $\bold c$ and $\bold d$ are arbitrary two vectors, we can further
simplify the obtained equations \thetag{6.15}, bringing them to the form
$$
\nabla_rb_s+b_r\,b'_s=\nabla_sb_r+b_s\,b'_r.\hskip -3em
\tag6.16
$$\par
    Next step consists in reducing the equations \thetag{5.2}. In order
to do it we substitute \thetag{5.5} into \thetag{5.2}. From \thetag{6.1}
we derive
$$
\sum^n_{s=1}\nabla_sA\,P^s_k=\sum^n_{s=1}\nabla_sa\,P^s_k+
\sum^n_{s=1}\sum^n_{r=1}\nabla_sb_r\,v^r\,P^s_k.\hskip -3em
\tag6.17
$$
Then from \thetag{6.3} and \thetag{6.4} we obtain the following two
relationships:
$$
\gather
\sum^n_{r =1}P^{qr}\,\tilde\nabla_qA=\sum^n_{s=1}b_s\,P^{sr},
\hskip -3em
\tag6.18\\
|\bold v|\sum^n_{r=1}\sum^n_{s=1}P^{qr}\,\tilde
\nabla_r\tilde\nabla_sA\,P^s_k=\left(a'+\shave{\sum^n_{r=1}}b'_r\,v^r
\right)P^q_k.\hskip -3em
\tag6.19
\endgather
$$
In \thetag{6.18} we have free index $r$, and in \thetag{6.19} we have
free index $q$. Let's multiply these two equalities \thetag{6.18} and
\thetag{6.19} and do contract with respect to indices $r$ and $q$ upon
multiplying the resulting equality by $g_{rq}$. This yields 
$$
|\bold v|\sum^n_{q=1}\sum^n_{r=1}\sum^n_{s=1}P^{qr}\,\tilde\nabla_qA
\,\tilde\nabla_r\tilde\nabla_sA\,P^s_k=
\left(a'+\shave{\sum^n_{r=1}}b'_r\,v^r\right)\sum^n_{s=1}b_s\,P^s_k.
$$
One more relationship is obtained from \thetag{6.4} upon multiplying
by $N^r\,A\,P^s_q$ and upon contracting with respect to $r$ and $s$:
$$
-\sum^n_{r=1}\sum^n_{s=1}N^r\,A\,\tilde\nabla_r\tilde\nabla_sA\,P^s_k=
-\left(a+\shave{\sum^n_{r=1}}b_r\,v^r\right)\sum^n_{s=1}b'_s\,P^s_k.
$$
Let's multiply \thetag{6.2} by $N^r$ and $P^s_k$, then contract it
with respect to $r$ and $s$:
$$
-|\bold v|\sum^n_{r=1}\sum^n_{s=1}N^r\,\nabla_r\tilde\nabla_sA
\,P^s_k=-\sum^n_{r=1}\sum^n_{s=1}v^r\,\nabla_rb_s\,P^s_k.
$$
Now, in order to write the result of substituting \thetag{5.5} into
the equations \thetag{5.2}, we have to add \thetag{6.17} and three
above equalities:
$$
\sum^n_{s=1}\left(\nabla_sa+b_s\,a'-a\,b'_s+\shave{\sum^n_{r=1}}
v^r\bigl(\nabla_sb_r+b_s\,b'_r-\nabla_rb_s-b_r\,b'_s\bigr)\right)
P^s_k=0.
$$
Let's take into account \thetag{6.16}, this leads to vanishing the
whole expression under summation with respect to $r$. As a result
we obtain the following equation:
$$
\sum^n_{s=1}\left(\nabla_sa+b_s\,a'-a\,b'_s\right)P^s_k=0.
\hskip -3em
\tag6.20
$$
Equations \thetag{6.20} are analogous to the equations \thetag{6.9},
the operation of contraction with components of projector $\bold P$
can be omitted:
$$
\nabla_sa+b_s\,a'=a\,b'_s.\hskip -3em
\tag6.21
$$
Arguments used in deriving the equations \thetag{6.21} are similar to
those used in deriving \thetag{6.16} from \thetag{6.9}.
\proclaim{Theorem 6.2} Force field $\bold F$ given by scalar ansatz
\thetag{5.1} corresponds to some Newtonian dynamical system admitting
the normal shift on the Riemannian manifold $M$ if and only if scalar
field $A$ in ansatz \thetag{5.1} is defined by formula \thetag{5.5},
while extended fields $a$ and $\bold b$ in \thetag{5.5} are fiberwise
spherically symmetric and satisfying the equations \thetag{6.16} and
\thetag{6.21}.
\endproclaim
Note that the equations \thetag{6.16} and \thetag{6.21} coincide with
reduced normality equations \thetag{5.7} and \thetag{5.6} we were to
derive.
\head
7. Analysis of reduced equations.
\endhead
     With the aim of further study of the equations \thetag{5.6} and
\thetag{5.7} let's express covariant derivatives in them through
partial derivatives. In order to do it we use formula \thetag{6.10}
and take into account symmetry of connection components:
$$
\gather
\left(\frac{\partial}{\partial x^s}+b_s\,\frac{\partial}{\partial v}
\right)a=\left(a\,\frac{\partial}{\partial v}\right)b_s,\hskip -2em
\tag7.1\\
\left(\frac{\partial}{\partial x^s}+b_s\,\frac{\partial}{\partial v}
\right)b_r=\left(\frac{\partial}{\partial x^r}+b_r\,\frac{\partial}
{\partial v}\right)b_s.\hskip -2em
\tag7.2
\endgather
$$
The equations \thetag{7.2} form closed system of equations with respect to
components $b_1,\,\ldots,\,b_n$ of covector field $\bold b$. We can study
them separately. Let's consider the differential operators in these
equations:
$$
\bold L_i=\frac{\partial}{\partial x^i}+b_i\,\frac{\partial}{\partial v}
\text{, \ where \ }i=1,\,\ldots,\,n.\hskip -2em
\tag7.3
$$
Now by means of direct calculations we can check that the equations
\thetag{7.2} are exactly the conditions of permutability of operators
\thetag{7.3}:
$$
[\bold L_s,\,\bold L_r]=0.\hskip -2em
\tag7.4
$$
Let $\Bbb R^+=(0,\,+\infty)$ be positive semiaxis on real axis $\Bbb R$.
Operators \thetag{4.3} have natural interpretation as vector fields on
the direct product of manifolds $M\times\Bbb R^+$. Let's complement
$\bold L_1,\,\ldots,\,\bold L_n$ by one more vector field $\bold L_{n+1}$,
which possibly is not commutating with fields $\bold L_1,\,\ldots,\,\bold
L_n$, but which should complete $\bold L_1,\,\ldots,\,\bold L_n$ up to a
moving frame on the manifold $M\times\Bbb R^+$. Each field $\bold L_i$
has its own local one-parametric group of local diffeomorphisms (see
\cite{19}) with parameter $y_i$:
$$
\varphi_i(y^i):\,M\times\Bbb R^+\to M\times\Bbb R^+.\hskip -2em
\tag7.5
$$
Let's fix some point $p_0\in M\times\Bbb R^+$ and let's consider composition
of such diffeomorphisms applied to the point $p_0$:
$$
p(y^1,\ldots,y^n,w)=\varphi_1(y^1)\compos\,\ldots\,\compos\varphi_n(y^n)
\compos\varphi_{n+1}(w)(p_0).\hskip -2em
\tag7.6
$$
In left hand side of the equality \thetag{7.6} we have the point $p$
parameterized by real numbers $y^1,\,\ldots,\,y^n,\,w$. This is equivalent
to defining local coordinates on $M\times\Bbb R^+$ in some neighborhood
of the point $p_0$. Permutability of vector fields \thetag{7.4} implies
permutability of first $n$ maps \thetag{7.5} in the composition
\thetag{7.6}. For the vector fields $\bold L_1,\,\ldots,\,\bold L_n$ this
fact yields the following expressions:
$$
\align
&\bold L_1=\frac{\partial}{\partial y^1}=\sum^n_{i=1}\frac{\partial x^i}
{\partial y^1}\,\frac{\partial}{\partial x^i}+\frac{\partial v}
{\partial y^1}\,\frac{\partial}{\partial v},\hskip -2em\\
&.\ .\ .\ .\ .\ .\ .\ .\ .\ .\ .\ .\ .\ .\ .\ .\ .\ .\ .\ .\ .\ .\
.\ .\ .\ .\ \hskip -2em
\tag7.7\\
&\bold L_n=\frac{\partial}{\partial y^n}=\sum^n_{i=1}\frac{\partial x^i}
{\partial y^n}\,\frac{\partial}{\partial x^i}+\frac{\partial v}
{\partial y^n}\,\frac{\partial}{\partial v}.\hskip -2em
\endalign
$$
Let's compare formulas \thetag{7.7} and \thetag{7.3} for vector fields
$\bold L_1,\,\ldots,\,\bold L_n$. This yields
$$
\frac{\partial x^i}{\partial y^k}=\delta^i_k=
\cases 1&\text{for \ }i=k,\\
0&\text{for \ }i\neq k.
\endcases\hskip -2em
\tag7.8
$$
From the same comparison for the functions $b_k$ in variables $y^1,\,
\ldots,\,y^n,\,w$ we get
$$
b_k=\frac{\partial v}{\partial y^k}.\hskip -2em
\tag7.9
$$
The relationships \thetag{7.8} show that newly constructed local
coordinates $y^1,\,\ldots,\,y^n,\,w$ on $M\times\Bbb R^+$ and initial
local coordinates $x^1,\,\ldots,\,x^n,\,v$ on this manifold are related
by means of only one function $V(y^1,\ldots,y^n,w)$:
$$
\cases
x^1=y^1,\,\ldots,\,x^n=y^n,\\
\vspace{0.5ex}
v=V(y^1,\ldots,y^n,w).
\endcases\hskip -3em
\tag7.10
$$
Inverse relation is also given by the only one function $W(x^1,\ldots,x^n,
v)$:
$$
\cases
y^1=x^1,\,\ldots,\,y^n=x^n,\\
\vspace{0.5ex}
w=W(x^1,\ldots,x^n,v).
\endcases\hskip -3em
\tag7.11
$$
Functions $V(y^1,\ldots,y^n,w)$ and $W(x^1,\ldots,x^n,v)$ in \thetag{7.10}
and in \thetag{7.11} are bound by the  obvious relationships that express
the fact that the changes of variables \thetag{7.10} and \thetag{7.11} are
inverse to each other:
$$
\aligned
&V(x^1,\ldots,x^n,W(x^1,\ldots,x^n,v))=v,
\vspace{0.5ex}
&W(x^1,\ldots,x^n,V(x^1,\ldots,x^n,w))=w.
\endaligned\hskip -3em
\tag7.12
$$
Let's use the following natural notations for partial derivatives of the
first order:
$$
\aligned
&V_i(x^1,\ldots,x^n,w)=\frac{\partial V(x^1,\ldots,x^n,w)}{\partial x^i},\\
\vspace{1.5ex}
&V_w(x^1,\ldots,x^n,w)=\frac{\partial V(x^1,\ldots,x^n,w)}{\partial w}.
\endaligned\hskip -3em
\tag7.13
$$
Analogous notations will be used for partial derivatives of the second order:
$$
\aligned
&V_{ij}(x^1,\ldots,x^n,w)=\frac{\partial^2 V(x^1,\ldots,x^n,w)}
{\partial x^i\,\partial x^j},\\
\vspace{1.5ex}
&V_{iw}(x^1,\ldots,x^n,w)=\frac{\partial^2 V(x^1,\ldots,x^n,w)}
{\partial x^i\,\partial w}.
\endaligned\hskip -3em
\tag7.14
$$
Now we can rewrite formula \thetag{7.9} in initial local coordinates
$x^1,\,\ldots,\,x^n,\,v$:
$$
b_k=V_k(x^1,\ldots,x^n,W(x^1,\ldots,x^n,v)).\hskip -3em
\tag7.15
$$
\proclaim{Theorem 7.1} Functions $b_1,\,\ldots,\,b_n$ satisfy nonlinear
system of partial differential equations \thetag{7.2} if and only if
they are determined by some function $V(x^1,\ldots,x^n,w)$ with non-zero
derivative $\partial V/\partial w$ according to the formula \thetag{7.15}.
\endproclaim
\demo{Proof} In theorem~7.1 we have two propositions. Direct proposition
is already proved: each solution of the system of equations \thetag{7.2}
is given by \thetag{7.15}. Conversely, suppose that some function
$V(x^1,\ldots,x^n,w)$ with non-zero derivative $\partial V/\partial w$
is chosen. From $\partial V/\partial w\neq 0$, relying on the theorem on
implicit functions (see \cite{51}, \cite{52}), we derive the existence
of the function $W(x^1,\ldots,x^n,v)$ such that it is bound with
$V(x^1,\ldots,x^n,w)$ by the relationships \thetag{7.12}. Differentiating
these relationships and taking into account the notations \thetag{7.13},
we derive
$$
\align
&\frac{\partial W}{\partial v}=\frac{1}{V_w(x^1,\ldots,x^n,W(x^1,\ldots,
x^n,v))}\hskip -3em
\tag7.16\\
\vspace{2.5ex}
&\frac{\partial W}{\partial x^s}=-\frac{V_s(x^1,\ldots,x^n,W(x^1,\ldots,
x^n,v))}{V_w(x^1,\ldots,x^n,W(x^1,\ldots,x^n,v))}\hskip -3em
\tag7.17
\endalign
$$
Let's substitute $V$ and $W$ into \thetag{7.15} and calculate functions
$b_1,\,\ldots,\,b_n$. Then by means of direct calculations we check that
the functions obtained satisfy differential equations \thetag{7.2}. Indeed,
here we have 
$$
\gather
\aligned
\frac{\partial b_r}{\partial x^s}=V_{rs}&(x^1,\ldots,x^n,W(x^1,
\ldots,x^n,v))\,+\\
&+\,V_{rw}(x^1,\ldots,x^n,W(x^1,\ldots,x^n,v))\,\frac{\partial W}
{\partial x^s},
\endaligned\hskip -3em
\tag7.18\\
\aligned
b_s\,\frac{\partial b_r}{\partial v}=V_s&(x^1,\ldots,x^n,W(x^1,\ldots,
x^n,v))\,\times\\
&\times\,V_{rw}(x^1,\ldots,x^n,W(x^1,\ldots,x^n,v))\,\frac{\partial W}
{\partial v}.
\endaligned\hskip -3em
\tag7.19
\endgather
$$
Let's add the equalities \thetag{7.18} and \thetag{7.19} and let's
take into account formulas \thetag{7.16} and \thetag{7.17} for 
derivatives. This yields the equality
$$
\left(\frac{\partial}{\partial x^s}+b_s\,\frac{\partial}{\partial v}
\right)b_r=V_{rs}(x^1,\ldots,x^n,W(x^1,\ldots,x^n,v)).\hskip -3em
\tag7.20
$$
Due to first formula \thetag{7.14} we can transpose indices $r$ and $s$
in right hand side of \thetag{7.20}, i\.~e\. $V_{rs}=V_{sr}$. This provides
the equations \thetag{7.2} for the functions \thetag{7.15} we have
constructed above.
\qed\enddemo
     Having constructed general solution for the equations \thetag{7.2},
now let's study the equations \thetag{7.1}. Let's implement the change of
variables \thetag{7.10} and transfer to variables $y^1,\,\ldots,\,y^n,\,w$.
In left hand side of \thetag{7.1} we have the same differential operator
$\bold L_s$ as in \thetag{7.2}. In variables $y^1,\,\ldots,\,y^n,\,w$ this
operator is written as $\bold L_s=\partial/\partial y^s$ (see relationships
\thetag{7.7}). Let's transform the operator in right hand side of
\thetag{7.1} to the variables $y^1,\,\ldots,\,y^n,\,w$:
$$
\frac{\partial}{\partial v}=\sum^n_{i=1}\frac{\partial y^i}
{\partial v}\,\frac{\partial}{\partial y^i}+
\frac{\partial w}{\partial v}\,\frac{\partial}{\partial w}=
\frac{\partial W}{\partial v}\,\frac{\partial}{\partial w}.
$$
For further transformation of the above expression for operator
$\partial/\partial v$ we use formula \thetag{7.16}. As a result we
get the following relationship:
$$
\frac{\partial}{\partial v}=\frac{1}{V_w(y^1,\ldots,y^n,w)}
\,\frac{\partial}{\partial w}.
$$
Now in variables $y^1,\,\ldots,\,y^n,\,w$ the equations \thetag{7.1} are
written as
$$
\frac{\partial a}{\partial y^s}=\frac{V_{sw}}{V_w}\,a.\hskip -3em
\tag7.21
$$
Here we used the above notations \thetag{7.13} and \thetag{7.14}. The
equations \thetag{7.21} can be easily solved if we rewrite them as follows:
$$
\frac{\partial}{\partial y^s}\left(\frac{a}{V_w}\right)=0.\hskip -3em
\tag7.22
$$
General solution of the equations \thetag{7.22} contains an arbitrary
function of one variable $h(w)$. It is given by the following formula:
$$
a=h(w)\ V_w(y^1,\ldots,y^n,w).\hskip -3em
\tag7.23
$$
Upon coming back to initial variables $x^1,\,\ldots,\,x^n,\,v$ from
\thetag{7.23} we obtain
$$
\aligned
a=h(&W(x^1,\ldots,x^n,v))\,\times\\
\vspace{1ex}
&\times\,V_w(x^1,\ldots,x^n,W(x^1,\ldots,x^n,v)).
\endaligned\hskip -3em
\tag7.24
$$\par
     With the aim of additional verification we substitute the
above expressions \thetag{7.24} and \thetag{7.15} into the equations
\thetag{7.1}. Let's do the appropriate calculations:
$$
\align
&\frac{\partial a}{\partial x^s}=\bigl(h'(W)\,V_w+h(W)\,V_{ww}\bigr)
\,\frac{\partial W}{\partial x^s}+h(W)\,V_{sw},\\
\vspace{1.5ex}
&b_s\,\frac{\partial a}{\partial v}=V_s\,\bigl(h'(W)\,V_w+h(W)\,V_{ww}
\bigr)\,\frac{\partial W}{\partial v}.
\endalign
$$
We add two above equalities and take into account formulas \thetag{7.16}
and \thetag{7.17} for partial derivatives $\partial W/\partial x^s$ and
$\partial W/\partial v$:
$$
\left(\frac{\partial}{\partial x^s}+b_s\,\frac{\partial}{\partial v}
\right)a=h(W)\,V_{sw}.\hskip -3em
\tag7.25
$$
Similar calculations for the right hand side of the equations \thetag{7.1}
yield
$$
\left(a\,\frac{\partial}{\partial v}\right)b_s=h(W)\,V_w\,
V_{sw}\frac{\partial W}{\partial v}.\hskip -3em
\tag7.26
$$
Comparing \thetag{7.25} with \thetag{7.26} and taking into account formulas
\thetag{7.16} for the derivative $\partial W/\partial v$ completes the proof
of the following theorem.
\proclaim{Theorem 7.2} Functions $b_1,\,\ldots,\,b_n$ and $a$ satisfy
nonlinear differential equations \thetag{7.1} and \thetag{7.2} if and
only if they are determined by formulas \thetag{7.15} and \thetag{7.24}.
\endproclaim
\head
8. General formula for force field.
\endhead
    Analyzing reduced normality equations above, we have found their
general solution. Now we are able to write {\bf formula} for the force
field of arbitrary Newtonian dynamical system {\bf admitting the normal
shift} on Riemannian manifold of the {\bf dimension $n\geqslant 3$}.
Let's substitute \thetag{5.5} into the formula for the components of
force field $\bold F$. Thereby we take into account that fields $a$ and
$\bold b$ determining scalar parameter $A$ are fiberwise spherically
symmetric:
$$
F_k=a\,N_k+|\bold v|\sum^n_{i=1}b_i\,\bigl(2\,N^i\,N_k-\delta^i_k\bigr).
\hskip -2em
\tag8.1
$$
Formula \thetag{8.1} completely determines the dependence of force
field on velocity vector $\bold v=v\,\bold N$. Taking into account
\thetag{7.15} and \thetag{7.24} we get
$$
\aligned
F_k&=h(W)\,\,V_w(x^1,\ldots,x^n,W)\,\,N_k\,+\\
\vspace{0.5ex}
&+\,|\bold v|\sum^n_{i=1}V_i(x^1,\ldots,x^n,W)\,\,\bigl(2\,N^i\,N_k
-\delta^i_k\bigr).
\endaligned\hskip -2em
\tag8.2
$$
Here $h(W)$ is an arbitrary function of one variable; through $V_i$ and
$V_w$ we denoted derivatives \thetag{7.13}, where functions $V(x^1,\ldots,
x^n,W)$ and $W(x^1,\ldots,x^n,v)$ are bound with each other by the
relationships \thetag{7.12}.\par
     Formula \thetag{8.2} for the force field of dynamical system admitting
the normal shift contains the arbitrariness determined by one function of
$(n+1)$ variables $V(x^1,\ldots,x^n,w)$. Arbitrariness determined by the
function $h(W)$ can be eliminated by means of gauge transformation that
changes $a$ but doesn't change $\bold b$:
$$
\align
&V(x^1,\ldots,x^n,w)\longrightarrow V(x^1,\ldots,x^n,\rho^{-1}(w)),\hskip 
-2em\\
\vspace{1ex}
&W(x^1,\ldots,x^n,v)\longrightarrow \rho(W(x^1,\ldots,x^n,v)),\hskip -2em
\tag8.3\\
\vspace{1ex}
&h(w)\longrightarrow h(\rho^{-1}(w))\,\,\rho'(\rho^{-1}(w)).
\endalign
$$
Hence transformation \thetag{8.3} doesn't change components of force
field $\bold F$, though it changes parameters $V$, $W$, and $h$ in
\thetag{8.2}. If $h(w)\neq 0$, we can choose function $\rho(w)$ such
that $h(w)\,\rho'(w)=1$. Upon doing gauge transformation \thetag{8.3}
in this case we obtain the following formula for force field $\bold F$:
$$
\aligned
F_k&=V_w(x^1,\ldots,x^n,W)\,\,N_k\,+\\
\vspace{0.5ex}
&+\,|\bold v|\sum^n_{i=1}V_i(x^1,\ldots,x^n,W)\,\,\bigl(2\,N^i\,N_k
-\delta^i_k\bigr).
\endaligned\hskip -2em
\tag8.4
$$
Formula \thetag{8.4} is almost as universal as formula \thetag{8.2}.
The only exception is the case $h=0$, which is not embraced by formula
\thetag{8.4}.
\head
9. Effectivization of general formula.
\endhead
     Formulas \thetag{8.2} and \thetag{8.4} determine components of the
force field of Newtonian dynamical system admitting the normal shift on
Riemannian manifold $M$. However, both these formulas have common fault.
They are ineffective since we are to use function $W(x^1,\ldots,x^n,v)$
determined implicitly by the equations \thetag{7.12}. With the aim to
get more effective formula we use the relationships \thetag{7.16} and
\thetag{7.17}. Let's rewrite these relationships as follows:
$$
\align
&V_w(x^1,\ldots,x^n,W(x^1,\ldots,x^n,v))=\frac{1}{\partial W
/\partial v},\hskip -2em
\tag9.1\\
\vspace{1.5ex}
&V_k(x^1,\ldots,x^n,W(x^1,\ldots,x^n,v))=-\frac{\partial W
/\partial x^k}{\partial W/\partial v}.\hskip -2em
\tag9.2
\endalign
$$
Substituting \thetag{9.2} into the formula \thetag{7.15}, for the
components of covector $\bold b$ we get
$$
b_k=-\frac{\partial W(x^1,\ldots,x^n,v)/\partial x^k}
{\partial W(x^1,\ldots,x^n,v)/\partial v}.\hskip -2em
\tag9.3
$$
Here $W(x^1,\ldots,x^n,v)$ can be understood as absolutely arbitrary
function provided the derivative in denominator of the fraction in
\thetag{9.3} is non-zero. The same function, upon substituting
\thetag{9.1} into \thetag{7.24}, determines the field $a$:
$$
a=\frac{h(W(x^1,\ldots,x^n,v))}{\partial W(x^1,\ldots,x^n,v)
/\partial v}.\hskip -2em
\tag9.4
$$
Let's substitute the expressions \thetag{9.3} and \thetag{9.4} into
the formula \thetag{8.1} for $\bold F$:
$$
\aligned
F_k&=\frac{h(W(x^1,\ldots,x^n,v))}{\partial W(x^1,\ldots,x^n,v)
/\partial v}\,\,N_k\,+\\
\vspace{1ex}
&-\,|\bold v|\sum^n_{i=1}\frac{\partial W(x^1,\ldots,x^n,v)/\partial x^i}
{\partial W(x^1,\ldots,x^n,v)/\partial v}\,\,\bigl(2\,N^i\,N_k
-\delta^i_k\bigr).
\endaligned\hskip -2em
\tag9.5
$$
Let's rewrite formula \thetag{9.5} in terms of covariant derivatives.
We formulate the result in form of the theorem.
\proclaim{Theorem 9.1} Newtonian dynamical system on Riemannian manifold
$M$ of the dimension $n\geqslant 3$ admits the normal shift if and only
if its force field $\bold F$ has the components determined by formula
$$
F_k=\frac{h(W)\,N_k}{W_v}-|\bold v|\sum^n_{i=1}\frac{\nabla_iW}{W_v}
\,\bigl(2\,N^i\,N_k-\delta^i_k\bigr),\hskip -2em
\tag9.6
$$
where $W$ is fiberwise spherically symmetric scalar field from extended
algebra of tensor fields on $M$ with non-zero derivative
$$
W_v=\frac{\partial W}{\partial v}=\sum^n_{i=1}N^i\,\,\tilde\nabla_i W
\neq 0,
$$
and $h=h(W)$ is an arbitrary function of one variable.
\endproclaim
\head
10. Kinematics of normal shift.
\endhead
     Having explicit formula \thetag{6.6} for the force field of the
dynamical admitting the normal shift, we are able to describe in details
the process of normal shift of a given hypersurface  $S$ along
trajectories of this dynamical system. According to the definition~3.2
one chooses some point $p_0$ on $S$, then on some part $S'=O_{\!S}(p_0)$
of hypersurface $S$ one should define the function $\nu$ that determines
the modulus of initial velocity on the trajectories of shift. At the point
$p_0$ this function $\nu$ should be normalized by the condition 
$$
\nu(p_0)=\nu_0,\hskip -3em
\tag10.1
$$
where $\nu_0$ is some nonzero number (see condition \thetag{3.3} above).
Let's choose local coordinates $x^1,\,\ldots,\,x^n$ on $M$ in some
neighborhood of the point $p_0$ and local coordinates $u^1,\,\ldots,\,
u^{n-1}$ on $S$ in a neighborhood of the same point. The choice of local
coordinates $u^1,\,\ldots,\,u^{n-1}$ determines coordinate tangent vectors
$\boldsymbol\tau_1,\,\ldots,\,\boldsymbol\tau_{n-1}$ to the hypersurface
$S$. Let $\tau^i_k$ be components of vector $\boldsymbol\tau_k$ in
coordinates $x^1,\,\ldots,\,x^n$ on $M$. If hypersurface $S$ is
defined parametrically
$$
\align
x^1&=x^1(u^1,\ldots,u^{n-1}),\hskip -3em\\
\vspace{-2pt}
.\ .\ &.\ .\ .\ .\ .\ .\ .\ .\ .\ .\ .\ .\ .\ .\ .
\hskip -3em\tag10.2\\
x^n&=x^n(u^1,\ldots,u^{n-1}),\hskip -3em
\endalign
$$
then components of coordinate tangent vectors $\boldsymbol\tau_1,\,\ldots,
\,\boldsymbol\tau_{n-1}$ are determined by derivatives of the functions
$x^i(u^1,\ldots,u^{n-1})$ in \thetag{10.2}:
$$
\tau^i_k=\frac{\partial x^i}{\partial u^k}.\hskip -3em
\tag10.3
$$
Function $\nu=\nu(p)=\nu(u^1,\ldots,u^{n-1})$ determines the initial
velocity on the trajectories of shift starting from $S$ (see formula
\thetag{2.3} above):
$$
\bold v(t)\,\hbox{\vrule height 8pt depth 8pt width 0.5pt}_{\,t=0}=
\nu(p)\cdot\bold n(p).\hskip -3em
\tag10.4
$$
Normal vector $\bold n(p)$ is determined up to a sign: $\bold n(p)
\to\pm\bold n(p)$. Therefore without loss of generality we can
assume that constant $\nu_0$ in \thetag{10.1} is positive. Then
function $\nu(p)$ in \thetag{10.4} is also positive. This means that
$$
|\bold v(t)|\,\hbox{\vrule height 8pt depth 10pt width 0.5pt}_{\,t=0}
=\nu(p).\hskip -3em
\tag10.5
$$\par
     The condition \thetag{10.4} provides normality of shift at initial
instant of time $t=0$. How to provide normality condition for other
instants of time $t\neq 0$\,? For this purpose we consider the solution
of Cauchy problem \thetag{2.1} for the system of differential equations
\thetag{1.2}, which describes Newtonian dynamical system with force
field $\bold F$. This solution is given by the functions \thetag{2.2}.
Let's write them as 
$$
\align
x^1&=x^1(u^1,\ldots,u^{n-1},t),\hskip -3em\\
\vspace{-2pt}
.\ .\ &.\ .\ .\ .\ .\ .\ .\ .\ .\ .\ .\ .\ .\ .\ .\ .\
\hskip -3em\tag10.6\\
x^n&=x^n(u^1,\ldots,u^{n-1},t).\hskip -3em
\endalign
$$
Functions \thetag{10.6} describe the shift $f_t\!:S\to S_t$. In
sufficiently small neighborhood $S'=O_{\!S}(p_0)$ of the point $p_0$
and for sufficiently small values of $t$ the map $f_t$ is  a
diffeomorphism: $f_t\!:S'\to S'_t$. Due to this diffeomorphism we can
carry local coordinates $u^1,\,\ldots,\,u^{n-1}$ from $S$ to $S_t$.
Then for any fixed $t$ the functions \thetag{10.6} can be treated as
parametric equations of hypersurface $S_t$ similar to the equations
\thetag{10.2}. Let's calculate derivatives \thetag{10.3} for the
functions \thetag{10.6}. Doing this, we define vectors $\boldsymbol\tau_1,
\,\ldots,\,\boldsymbol\tau_{n-1}$ tangent to all hypersurfaces $S_t$:
$$
\boldsymbol\tau_k=\boldsymbol\tau_k(u^1,\ldots,u^{n-1},t).\hskip -3em
\tag10.7
$$
Denote by $\varphi_i$ the scalar product of the vector \thetag{10.7}
and the vector of velocity:
$$
\varphi_k=(\boldsymbol\tau_k\,|\,\bold v).\hskip -3em
\tag10.8
$$
In thesis \cite{17} the functions $\varphi_1,\,\ldots,\,\varphi_{n-1}$
determined by formula \thetag{10.8} were called {\bf functions of
deviation}. Such functions play an important role in deriving the
normality equations \thetag{3.1} and \thetag{3.2}, since the condition
of normality for the shift $f_t\!:S\to S_t$ consists in identical
vanishing of all functions \thetag{10.8}:
$$
\varphi_k=\varphi_k(u^1,\ldots,u^{n-1},t)=0.\hskip -3em
\tag10.9
$$
From identical vanishing of the functions of deviation \thetag{10.9}
it follows that their time derivatives at the initial instant of
time $t=0$ are zero:
$$
\dot\varphi_k\,\hbox{\vrule height 8pt depth 8pt width 0.5pt}_{\,t=0}
=0.\hskip -3em
\tag10.10
$$
Moreover from the relationships \thetag{10.9} it follows that $\varphi_k$
for $t=0$ are zero as well:
$$
\varphi_k\,\hbox{\vrule height 8pt depth 8pt width 0.5pt}_{\,t=0}
=0.\hskip -3em
\tag10.11
$$
In papers \cite{6} and \cite{7} it was shown that for Newtonian dynamical
systems admitting the normal shift (in the sense of definition~3.2 and
theorem~3.1) the conditions \thetag{10.10} and \thetag{10.11} are not
only necessary, but also sufficient for identical vanishing of all
functions of deviation \thetag{10.8}.\par
     The conditions \thetag{10.11} are trivial consequences of
\thetag{10.4}. They give no information on how to choose the function 
$\nu(p)$ on $S$. Therefore let's consider the conditions \thetag{10.10}.
Let's calculate $\dot\varphi$ by differentiating \thetag{10.8}:
$$
\dot\varphi_k=\nabla_t\varphi_k=(\nabla_t\boldsymbol\tau_k\,|\,\bold v)+
(\boldsymbol\tau_k\,|\,\bold F).\hskip -3em
\tag10.12
$$
In deriving \thetag{10.12} we took into account the equations of
dynamics \thetag{1.2} written as $\nabla_t\bold v=\bold F$. Now let's
calculate the derivative $\nabla_t\boldsymbol\tau_k$. It's the vector
with components
$$
\nabla_t\tau^i_k=\frac{\partial\tau^i_k}{\partial t}+
\sum^n_{r=1}\sum^n_{s=1}\Gamma^i_{rs}\,v^r\,\tau^s_k.\hskip -3em
\tag10.13
$$
If we take into account formulas \thetag{10.3} determining $\tau^i_k$,
then from \thetag{10.13} we obtain
$$
\nabla_t\boldsymbol\tau^i_k=\frac{\partial v^i}{\partial u^k}+
\sum^n_{r=1}\sum^n_{s=1}\Gamma^i_{rs}\,v^r\,\frac{\partial x^s}
{\partial u^k}=\nabla_{u^k}v^i.\hskip -3em
\tag10.14
$$
In vectorial form \thetag{10.14} is written as $\nabla_t\boldsymbol
\tau_k=\nabla_{u^k}\bold v$. That is $\nabla_t\boldsymbol\tau_k$
coincides with covariant derivative of vector function $\bold v(u^1,
\ldots,u^{n-1},t)$ with respect to parameter $u^k$ along $k$-th
coordinate line on hypersurface $S_t$. Let's substitute the obtained
expression for $\nabla_t\boldsymbol\tau_k$ into the formula
\thetag{10.12}. This yields
$$
\dot\varphi_k=(\nabla_{u^k}\bold v\,|\,\bold v)+
(\boldsymbol\tau_k\,|\,\bold F).\hskip -3em
\tag10.15
$$
Further, we take into account the obvious relationship $\nabla_{u^k}
(\bold v\,|\,\bold v)=2\,(\nabla_{u^k}\bold v\,|\,\bold v)$ and formula
\thetag{10.5}, which determines modulus of velocity vector for $t=0$.
Then we can bring formula \thetag{10.15} for the derivative
$\dot\varphi_k$ to the following form:
$$
\dot\varphi_k\,\hbox{\vrule height 8pt depth 8pt width 0.5pt}_{\,t=0}=
\nu\,\frac{\partial\nu}{\partial u^k}+(\boldsymbol\tau_k\,|\,\bold F).
\hskip -3em
\tag10.16
$$
Now, due to \thetag{10.16}, the relationships \thetag{10.10} \pagebreak
for derivatives $\dot\varphi_k$ are written as partial differential
equations for the function $\nu=\nu(u^1,\ldots,u^{n-1})$ on $S$:
$$
\frac{\partial\nu}{\partial u^k}=-\nu^{-1}\,(\bold F\,|\,
\boldsymbol\tau_k).\hskip -3em
\tag10.17
$$
Let's substitute force field \thetag{9.6} into the equations \thetag{10.17}
and take into account the fact that vector $\bold N$ for $t=0$ coincides
with unitary normal vector $\bold n(p)$ on $S$. Upon rather simple
calculations this yields
$$
\frac{\partial\nu}{\partial u^k}=-\sum^n_{i=1}\frac{\nabla_iW}{W_v}
\,\tau^i_k.\hskip -3em
\tag10.18
$$
Let's multiply the equation \thetag{10.18} by $W_v$ and transfer the
sum from left to right hand side of this equation. Moreover, let's
write explicitly all derivatives:
$$
\frac{\partial W}{\partial v}\,\frac{\partial\nu}{\partial u^k}
+\sum^n_{i=1}\frac{\partial W}{\partial x^i}\,\frac{\partial x^i}
{\partial u^k}=0.\hskip -3em
\tag10.19
$$
It's not difficult to see that left hand side of \thetag{10.19} is the
derivative of the function $W(x^1,\ldots,x^n,v)$ with respect to $u^k$
upon substituting the functions \thetag{10.2} for $x^1,\,\ldots,\,x^n$
and the function $\nu(u^1,\ldots,u^{n-1})$ for $v$. Therefore the equations
\thetag{10.19} are easily integrated in form of functional equation
$$
W(x^1(p),\ldots,x^n(p),\nu(p))=W_0=\const,\hskip -3em
\tag10.20
$$
which determines the function $\nu=\nu(p)=\nu(u^1,\ldots,u^{n-1})$ in
implicit form. The value of constant $W_0$ in \thetag{10.20} is fixed
by normalizing condition \thetag{10.1}:
$$
W_0=W(x^1(p_0),\ldots,x^n(p_0),\nu_0).\hskip -3em
$$
\proclaim{Theorem 10.1} In order to construct the normal shift of
hypersurface $S$, given in parametric form by functions $x^1(p),\,\ldots,
\,x^n(p)$ from \thetag{10.2}, along trajectories of Newtonian dynamical
system with force field \thetag{9.6} one should determine the function
$\nu(p)$ in \thetag{10.4} by means of the equation \thetag{10.20}.
\endproclaim
     Having constructed the normal shift $f_t\!:\,S\to S_t$ along
trajectories of dynamical system with force field \thetag{9.6}, we
obtain the family of hypersurfaces $S_t$. By changing the initial
instant for counting the time $t\to t+t_0$ we can treat each hypersurface
of this family as initial hypersurface. Therefore on each of them
the following equality similar to \thetag{10.20} is fulfilled:
$$
W(p,|\bold v|)=W_0(t)=\const.\hskip -3em
\tag10.21
$$
Note that the values of constants $W_0(t)$ in \thetag{10.21} can be
different on different hypersurfaces $S_t$. Let's calculate the dynamics
of $W_0(t)$ in $t$. First find the time dynamics of the modulus of
velocity vector on the trajectories of shift:
$$
\frac{d\,|\bold v|}{dt}=\nabla_tv=
\frac{(\bold v\,|\,\nabla_t\bold v)}{v}=(\bold N\,|\,\bold F)
=\sum^n_{k=1}N^k\,F_k.\hskip -3em
\tag10.22
$$
Then substitute \thetag{9.6} into \thetag{10.22}. As a result of
this substitution we obtain
$$
\frac{dv}{dt}=\frac{h(W)}{W_v}-\sum^n_{i=1}\frac{\nabla_iW}
{W_v}\,\frac{dx^i}{dt}.\hskip -3em
\tag10.23
$$
Now let's multiply \thetag{10.23} by $W_v$ and transfer the
sum from left to right hand side of this equation. Moreover, let's
write explicitly all derivatives:
$$
\sum^n_{i=1}\frac{\partial W}{\partial x^i}\,\frac{dx^i}{dt}+
\frac{\partial W}{\partial v}\,\frac{dv}{dt}=h(W).\hskip -3em
\tag10.24
$$
In left hand side of \thetag{10.24} we see the time derivative
of the scalar field on the trajectories of normal shift. Hence
\thetag{10.24} is the required equation that determines time
dynamics of constants $W_0(t)$ in \thetag{10.21}. Let's
write this equation as follows:
$$
\frac{dW_0}{dt}=h(W_0).\hskip -3em
\tag10.25
$$
If function $h(w)$ in \thetag{9.6} is zero, then due to \thetag{10.25}
the field $W$ in \thetag{10.21} not only is constant on each separate
hypersurface $S_t$, but has equal values on all such hypersurfaces.
For $h(w)\neq 0$ differential equation \thetag{10.25} is easily integrated.
So, knowing the value of $W$ on $S$, we can find its value on any one
of hypersurfaces $S_t$.\par
\head
11. Coordinates associated with Newtonian normal shift.
\endhead
    It is known that the construction of geodesic normal shift of
hypersurface $S$ in Riemannian manifold $M$ provides some special
choice of local coordinates in a neighborhood of $S$. They are
called {\bf semigeodesic coordinates} (see \cite{20} or \cite{53}).
Newtonian normal shift of hypersurface $S$ also can provide some
special choice of local coordinates in $M$. Let $f_t\!:S\to S_t$ be
the normal shift of $S$ along trajectories of Newtonian dynamical
system with force field \thetag{9.6}, function $\nu(p)$ for which
is fixed by normalizing condition \thetag{10.1} at the point $p_0
\in S$ (without loss of generality we can assume that $\nu_0>0$).
Then some neighborhood of the point $p_0$ foliates into the union
of not intersecting parts of hypersurfaces $S_t$. Choosing local
coordinates $u^1,\,\ldots,\,u^{n-1}$ on $S$ in a neighborhood of
$p_0$, we can carry them from $S$ to $S_t$ by means of shift
diffeomorphism $f_t\!:S\to S_t$. Therefore the set of $n$ quantities
$u^1,\,\ldots,\,u^{n-1}$, and $t$ can be considered as local coordinates
in $M$ in some neighborhood of the point $p_0$. Such coordinates are
called {\bf associated with normal shift} $f_t\!:S\to S_t$. If we denote
associated coordinates by $x^1,\,\ldots,\,x^n$, then the functions
$x^i(u^1,\ldots,u^{n-1},t)$ in \thetag{10.6} become extremely simple:
$$
\align
&x^1(u^1,\ldots,u^{n-1},t)=u^1,\hskip -3em\\
\vspace{-2pt}
&.\ .\ .\ .\ .\ .\ .\ .\ .\ .\ .\ .\ .\ .\ .\ .\ .\ .\
\hskip -3em\tag11.1\\
&x^{n-1}(u^1,\ldots,u^{n-1},t)=u^{n-1},\hskip -3em\\
&x^n(u^1,\ldots,u^{n-1},t)=t.
\endalign
$$
Differentiating functions \thetag{11.1} according to \thetag{10.3},
\pagebreak we obtain the components of vectors $\boldsymbol\tau_1,\,
\ldots,\,\boldsymbol\tau_{n-1}$ tangent to the hypersurfaces $S_t$
in associated coordinates:
$$
\tau^i_k=\delta^i_k=\cases 1&\text{for \ }i=k,\\
0&\text{for \ }i\neq k.\endcases
$$
Differentiating functions \thetag{11.1} in $t$, we find components of
velocity vector $\bold v$:
$$
v^i=\delta^i_n=\cases 1&\text{for \ }i=n,\\
0&\text{for \ }i\neq n.\endcases\hskip -3em
\tag11.2
$$
Denote by $\nu$ the modulus of velocity vector on trajectories of
shift:
$$
\nu=\nu(u^1,\ldots,u^{n-1},t)=|\bold v|.\hskip -3em
\tag11.3
$$
Doing this, we extend the domain of function $\nu$ from \thetag{10.5},
which is initially defined only for $t=0$ on initial hypersurface $S$.
From \thetag{11.2} we obtain the relationship
$$
\nu=|\bold v|=\sqrt{g_{nn}(u^1,\ldots,u^{n-1},t)}.\hskip -3em
\tag11.4
$$
On the other hand, the modulus of velocity vector can be calculated
from the functional equation \thetag{10.21}. One can solve this
equation in explicit form by using the function $V(x^1,\ldots,x^n,w)$
from \thetag{7.10}:
$$
|\bold v|=V(u^1,\ldots,u^{n-1},t,W_0(t)).\hskip -3em
\tag11.5
$$
Comparing \thetag{11.4} and \thetag{11.5}, we get the following formula
for diagonal component $g_{nn}$ of metric tensor in local coordinates
associated with normal shift $f_t\!:S\to S_t$:
$$
g_{nn}={V(x^1,\ldots,x^n,W_0(x^n))}^2.\hskip -3em
\tag11.6
$$
Following non-diagonal components of $g_{ij}$ are zero due to normality
of shift:
$$
g_{nk}=0\text{\ \ for all \ } k=1,\,\ldots,\,n-1.\hskip -3em
\tag11.7
$$
Function of one variable $W_0(t)$ from \thetag{11.5} and \thetag{11.6}
is determined as the solution of ordinary differential equation
\thetag{10.25} fixed by initial condition 
$$
W_0(t)\,\hbox{\vrule height 8pt depth 8pt width 0.5pt}_{\,t=0}=
W(p_0,\nu_0).\hskip -3em
\tag11.8
$$
Initial condition \thetag{11.8} follows from \thetag{10.21}, from
the relationship \thetag{10.5}, and from normalizing condition
\thetag{10.1}.\par
    The relationships \thetag{11.6} and \thetag{11.7} mean that
matrix formed by components of metric tensor in local coordinates
$x^1,\,\ldots,\,x^n$ associated with Newtonian normal shift
$f_t\!:S\to S_t$ has blockwise-diagonal structure:
$$
g=\Vmatrix
\lower 12pt\hbox{\vbox{\hrule\hbox to 40pt{\vrule height 18pt depth 13pt\hss 
$G$\hss
\vrule height 18pt depth 13pt}\hrule}} &
\matrix 0\\ \vdots\\ 0\endmatrix\\
\vspace{5pt}
\matrix 0 &\hdots &0\endmatrix & g_{nn}
\endVmatrix.\hskip -3em
\tag11.9
$$
Trajectories of normal shift $f_t\!:S\to S_t$ correspond to the variation
of parameter $t$ in \thetag{11.1} by fixed values of parameters $u^1,\,
\ldots,\,u^{n-1}$. Let's write the equation of Newtonian dynamics of points
$\nabla_t\bold v=\bold F$ for such trajectories. In associated local
coordinates this vectorial equation reduces to the series of scalar
equations. Taking into account the relationship \thetag{11.2}, we obtain
$$
\Gamma^k_{nn}=F^k,\qquad k=1,\,\ldots,\,n.\hskip -3em
\tag11.10
$$
Let's lower the upper index $k$ in \thetag{11.10}. Thereby we take
into account blockwise-diagonal structure of matrix of metric tensor
\thetag{11.9} and explicit formula \thetag{6.14} for components of
metric connection. This yields
$$
\align
\frac{\partial g_{nn}}{\partial x^k}&=-2\,F_k\text{\ \ for \ }k<n,
\hskip -3em
\tag11.11\\
\vspace{2ex}
\frac{\partial g_{nn}}{\partial x^n}&=2\,F_n\text{\ \ \ \ for \ }k=n.
\hskip -3em
\tag11.12
\endalign
$$
The quantity $g_{nn}$ in left hand sides of \thetag{11.11} and
\thetag{11.12} is determined by formula \thetag{11.6}. Covariant
components of force vector are determined by formula \thetag{9.6},
or equivalent formula \thetag{8.2}. By substituting \thetag{11.6}
and \thetag{8.2} into \thetag{11.11} and into \thetag{11.12} we
take into account that components of velocity vector on trajectories
of shift  are determined by formulas \thetag{11.2}, while its modulus
is determined by formula \thetag{11.3}. Therefore for components
of unitary vector $\bold N$ we have
$$
\xalignat 2
&N^k=\frac{\delta^k_n}{|\bold v|},&&N_k=|\bold v|\,\delta^n_k.
\endxalignat
$$
If we remember all circumstances listed above, then by substituting
\thetag{11.6} and \thetag{8.2} into the equalities \thetag{11.11} and
\thetag{11.12} we find that these equalities turn to identities. So
we get no restrictions for the choice of functions $h(w)$ and
$V(x^1,\ldots,x^n,w)$. This is not surprising, since all restrictions
due to normality of shift $f_t\!:S\to S_t$ are already handled by
normality equations \thetag{3.1} and \thetag{3.2}, and by explicit
formula \thetag{8.2} that follows from these equations. As for the
normalizing condition \thetag{3.3} for $\nu$, it is provided by the
relationship \thetag{11.6}, by the equation \thetag{10.25} for the
function $W_0(t)$, and by initial condition \thetag{11.8}.\par
     The above result gives the answer to one of the questions by
A.~V.~Bolsinov and A.~T.~Fomenko. It is formulated as follows: {\it
to what extent the normal shift of some particular hypersurface
$f_t\!:S\to S_t$ characterizes the structure of force field of
dynamical system admitting the normal shift}\,? The answer is:
{\bf yes, it characterizes, but partially; it doesn't determine
$\bold F$ completely}.\par
     Indeed, if normal shift $f_t\!:S\to S_t$ is already constructed,
then constructing associated coordinates $x^1,\,\ldots,\,x^n$ in a
neighborhood of $S$ reduces to the choice of local coordinates
$u^1,\,\ldots,\,u^{n-1}$ on initial hypersurface $S$. Diagonal component
$g_{nn}$ of metric tensor in these coordinates determines the function
of $n$ variables 
$$
V(x^1,\ldots,x^n,W_0(x^n))=\sqrt{g_{nn}(x^1,\ldots,x^n)}
\hskip -3em
\tag11.13
$$
(see the relationship \thetag{11.6} above). \pagebreak But by function
\thetag{11.3} one cannot reconstruct the function of $n+1$ variables
$V(x^1,\ldots,x^n,w)$, which is contained in the formula \thetag{8.2}
for components of force field $\bold F$.\par
     If we suppose that $h(w)=0$ in formula \thetag{8.2}, then from
\thetag{10.25} it follows that $W_0(x^n)=W_0=\const$. In this case
formula \thetag{11.13} simplifies to
$$
V(x^1,\ldots,x^n,W_0)=\sqrt{g_{nn}(x^1,\ldots,x^n)}.
\hskip -3em
\tag11.14
$$
But for the fixed normal shift $f_t\!:S\to S_t$ the constant $W_0$
is strictly fixed. Therefore by \thetag{11.14} we cannot reconstruct
the function $V(x^1,\ldots,x^n,w)$ in whole.\par
     Now suppose that $h(w)\neq 0$. Let's consider gauge transformations
\thetag{8.3} that do not change force field \thetag{8.2}. Transformations
\thetag{8.3} are supplemented by the rule for transforming $W_0(t)$. It
looks like
$$
W_0(t)\longrightarrow \rho(W_0(t)).\hskip -3em
\tag11.15
$$
Gauge transformations \thetag{8.3} supplemented by the additional rule
\thetag{11.15} preserve not only the force field $\bold F$ of dynamical
system, but the maps of normal shift $f_t$ as well. Therefore they do
not change the choice of associated local coordinates and the function
$g_{nn}$ in right hand side of \thetag{11.13}. Invariance of the left
hand side of \thetag{11.13} with respect to these transformations is
easily checked by direct calculations. At the expense of the gauge
transformations \thetag{8.3} the case when $h(w)\neq 0$ can be reduced
to the case $h(w)=1$ (see comment preceding formula \thetag{8.4}).
For $h(w)=1$ by means of integrating \thetag{10.25} we get $W_0(t)=W_0+t$,
where $W_0=\const$. This reduces \thetag{11.13} to the following form:
$$
V(x^1,\ldots,x^n,W_0+x^n)=\sqrt{g_{nn}(x^1,\ldots,x^n)}.
\hskip -3em
\tag11.16
$$
Constant $W_0$ in \thetag{11.6} is strictly fixed for the fixed normal
shift $f_t\!:S\to S_t$. Therefore by \thetag{11.16} one cannot reconstruct
the function \thetag{11.16} in whole.\par
     {\bf Consider a simple example}. Let $M=\Bbb R^n$ be euclidean 
space with standard metric, and let $S$ be hyperplane given by the
equation $x^n=0$. Then the relationships \thetag{11.1} determine
parametrization of $S$ and define the normal shift of this hypersurface
$f_t\!:S\to S_t$, being the parallel displacement of $S$ along $n$-th
coordinate axis. Thereby $g_{nn}=1$.\par
     1. The above shift can be implemented by Newtonian dynamical system
with identically zero force field $\bold F(\bold r,\bold v)=0$. This
corresponds to the choice $h(w)=0$ and $V(x^1,\ldots,x^n,w)=w$ in
formula \thetag{8.2}, and to the choice $W_0=1$ in formula \thetag{11.14}
respectively.\par
    2. Function $V(x^1,\ldots,x^n,w)$ can be changed, keeping $h(w)=0$
and $W_0=1$ meanwhile. Let's set $V(x^1,\ldots,x^n,w)=w+(1-w)\cdot
\varphi(x^1,\ldots,x^n)$. Then the above normal shift of hyperplane $S$
will be implemented by Newtonian dynamical system, force field of which
is non-zero:
$$
\bold F(\bold r,\bold v)=\frac{1-|\bold v|}{1-\varphi}\cdot\frac{2\,
(\bold v\,|\,\nabla\varphi)\cdot\bold v-|\bold v|^2\cdot\nabla\varphi}
{|\bold v|}.
$$
Here $\varphi=\varphi(\bold r)$ is an arbitrary function of coordinates
$x^1,\,\ldots,\,x^n$, for which is natural to assume, that its values are
distinct from $1$.\par
    3. Taking $h(w)=1$, we can choose $W(x^1,\ldots,x^n)=w-x^n$ in formula
\thetag{8.2}, and $W_0=1$ in formula \thetag{11.16} respectively. For
the force field $\bold F$ this yields 
$$
\bold F(\bold r,\bold v)=\frac{\bold v}{|\bold v|}-
\frac{2\,(\bold v\,|\,\bold M)\cdot\bold v-|\bold v|^2\cdot\bold M}
{|\bold v|}.
$$
The example, which we have just examined, confirms our conclusion that
knowing the normal shift $f_t\!:S\to S_t$ of some particular hypersurface
is not sufficient for to determine the force field of dynamical system
implementing this shift, even if we know that this system belong to the
class of systems admitting the normal shift.\par
\head
12. Blowing up the points. Generalization of theory as
suggested by A.~V.~Bolsinov and A.~T.~Fomenko.
\endhead
\parshape 22 0pt 360pt 0pt 360pt 0pt 360pt 0pt 360pt 0pt 360pt
0pt 360pt 0pt 360pt 160pt 200pt 
160pt 200pt 160pt 200pt 160pt 200pt 160pt 200pt 160pt 200pt 160pt 200pt
160pt 200pt 160pt 200pt 160pt 200pt 160pt 200pt 160pt 200pt
160pt 200pt 160pt 200pt 0pt 360pt
    Before now, studying normal shift $f_t\!:S\to S_t$, we restricted
ourselves to the case of smooth hypersurfaces and took parameter $t$
small enough for hypersurfaces $S_t$ to be non-singular as well.
However, one case with singularity appears to be interesting now. This
is the case when hypersurface $S_t$ collapses into a point at a time
for some $t=t_0$. By reverting the direction of time we can speak about
blowing up the point. Moreover, without loss of generality we can assume
that $t_0=0$. In this case we have singular initial hypersurface
$S=\{p_0\}$ consisting of only one point $p_0$, and a fan-shaped pencil
of trajectories coming out from this point (see figure~12.1). Velocity
vectors $\bold v=\bold v(0)$ on these trajectories corresponding to
the time instant $t=0$ belong to the tangent space $T_{p_0}(M)$. They
determine a hypersurface $s$ in the fiber of tangent bundle $TM$ over
the point $p_0\in M$. \vadjust{\vskip -65pt\hbox to 0pt{\kern 5pt
\hbox{\special{em:graph Eds-01c.gif}}\hss}\vskip 65pt}It can be
understood as ``limiting variety'' for hypersurfaces $S_t$ in
``infinitesimal scale'':
$$
s=\lim_{t\to 0}\frac{S_t}{t}.\hskip -2em
\tag12.1
$$
If initial values of velocity vectors on all trajectories at the point
$p_0$ are non-zero, then $s$ have topology of $(n-1)$-dimensional sphere.
The same topology is inherited by all hypersurfaces $S_t$ for sufficiently
small values of parameter $t$. Let $\sigma$ be the unit sphere in the fiber
of tangent bundle over the point $p_0$, let $q$ be a point of this sphere,
and let $\bold n(q)$ be radius-vector of the point $q$ in $T_{p_0}(M)$.
Then radius-vectors of the points on the hypersurface $s$ are given by
formula $\bold v(q)=\nu(q)\cdot\bold n(q)$, where $\nu=\nu(q)$ is some
positive function on unit sphere $\sigma$, while trajectories coming out
from the point $p_0$ are determined by initial data
$$
\xalignat 2
&\quad x^k\,\hbox{\vrule height 8pt depth 8pt width 0.5pt}_{\,t=0}
=x^k(p_0),
&&\dot x^k\,\hbox{\vrule height 8pt depth 8pt width 0.5pt}_{\,t=0}=
\nu(q)\cdot n^k(q)\hskip -3em
\tag12.2
\endxalignat
$$
for the equations of Newtonian dynamics \thetag{1.2}. Components of
force vector in \thetag{1.2} are determined by formulas \thetag{8.2}.
Due to normality of shift $f_t\!:S\to S_t$ all results obtained above
in sections~10 and 11 \pagebreak remain valid for non-singular
hypersurfaces $S_t$ with $t=0$. Formula \thetag{11.5} from section~11
now is written as follows:
$$
|\bold v|=V(x^1(t,q),\ldots,x^n(t,q),W_0(t)).\hskip -3em
\tag12.3
$$
We can return to initial form of formula \thetag{11.5} if we denote
by $u^1,\,\ldots,\,u^{n-1}$ local coordinates of the point $q$ on
unit sphere $\sigma$.\par
     Modulus of initial velocity $|\bold v|=\nu(q)$, which is contained
in formula \thetag{12.2}, can be found by passing to the limit $t\to 0$
in formula \thetag{12.3}:
$$
\nu(q)=V(x^1(p_0),\ldots,x^n(p_0),W_0).\hskip -3em
\tag12.4
$$
Here $x^1(p_0),\,\ldots,\,x^n(p_0)$ are coordinates of the point
$p_0$, and $W_0$ is initial value of function $W_0(t)$ for $t=0$.
The function $W_0(t)$ itself is determined as solution of ordinary
differential equation \thetag{10.25}. Note that right hand side of
\thetag{12.4} doesn't depend on $q$, i\.~e\. function $\nu(q)$ is
constant:
$$
\nu(q)=V(p_0,W_0)=\nu_0=\const.\hskip -3em
\tag12.5
$$
Constant $W_0$ can be expressed through $\nu_0$ if we take into
account \thetag{7.12}:
$$
W_0=W(p_0,\nu_0)=\const.\hskip -3em
\tag12.6
$$
On account of the relationships \thetag{12.5} and \thetag{12.6} we
can rewrite \thetag{12.2} as follows:
$$
\xalignat 2
&\quad x^k\,\hbox{\vrule height 8pt depth 8pt width 0.5pt}_{\,t=0}
=x^k(p_0),
&&\dot x^k\,\hbox{\vrule height 8pt depth 8pt width 0.5pt}_{\,t=0}=
\nu_0\cdot n^k(q).\hskip -3em
\tag12.6
\endxalignat
$$
Thus, we can state a theorem that follows from the results of
sections~10 and 11 by passing to the limit \thetag{12.1}.
\proclaim{Theorem 12.1} Suppose that on Riemannian manifold $M$ some
Newtonian dynamical system admitting the normal shift is defined,
i\.~e\. we have the system with force field $\bold F$ given by formula
\thetag{8.2}. Then for any point $p_0\in M$ and for any positive
constant $\nu_0>0$ initial data \thetag{12.6} determine 
the normal blow-up $f_t\!:p_0\to S_t$ of the point $p_0$ along
trajectories of this dynamical system.
\endproclaim
     Consideration of normal blow-ups for separate point of the
manifold $M$ gives the opportunity for further development of the
theory of dynamical systems admitting the normal shift. We can
formulate the following definition similar to definition~3.1.
\definition{Definition 12.1} Newtonian dynamical system \thetag{1.2}
with force field $\bold F$ on Riemannian manifold $M$ is called a system
{\bf admitting the normal blow-ups of points} if for any point $p_0\in M$
and for arbitrary positive constant $\nu_0>0$ initial data \thetag{12.6}
determine the normal blow-up $f_t\!:p_0\to S_t$ of the point $p_0$ along
trajectories of this dynamical system.
\enddefinition
    The idea of constructing new theory on the base of definition~12.1
was suggested by A.~V.~Bolsinov and A.~T.~Fomenko in February of 2000
in the seminar at Moscow State University during the discussion on the
results of thesis \cite{17}. Theorem~12.1 shows that dynamical systems
with force field \thetag{8.2} are included into the framework of new
theory. \pagebreak But, possibly, one can find some new dynamical systems,
which aren't embraced by formula \thetag{8.2}. We can compare this situation
with that of the theory of distributions, where narrowing class of test
functions extends the class of distributions. Here, narrowing
class of initial hypersurfaces in the construction of normal shift
$f_t\!:S\to S_t$ to singular one-point sets, we have a good chance to
extend class of dynamical systems that can implement such shift. Is it
really so\,? The answer to this question can be given only as a result
of constructing new theory. But this falls out of the limits of this
paper.
\head
13. On the problems of metrizability.
\endhead
    Problem of metrizability arose on initial stage of developing
theory of dynamical systems admitting the normal shift by testing
theory for non-triviality. It was noted (see papers \cite{5} and
\cite{7}) that if metric $\tilde\bold g$ is conformally equivalent
to the basic metric $\bold g$ of Riemannian manifold $M$, i\.~e\.
if we have 
$$
\tilde\bold g=e^{-2f}\,\bold g,\hskip -3em
\tag13.1
$$
then geodesic flow of metric $\tilde\bold g$ is a dynamical system
admitting the normal shift with respect to metric $\bold g$. Its
force field is given by formula 
$$
F_k=-|\bold v|^2\,\nabla_kf+2\sum^n_{s=1}
\nabla_sf\,v^s\,v_k.\hskip -3em
\tag13.2
$$
Here $\nabla f$ is a gradient of scalar field $f$ determining conformal
factor $e^{-2f}$ in \thetag{13.1}.
\definition{Definition 13.1} Newtonian dynamical system on Riemannian
manifold is called {\bf metrizable system} if it inherits trajectories
of the system with force field \thetag{13.2}.
\enddefinition
    Trajectory inheriting and trajectory equivalence for two dynamical
systems are understood in the sense of the following definitions.
\definition{Definition 13.2} Suppose that on the Riemannian manifold
$M$ two Newtonian dynamical systems are defined with force fields
$\bold F$ and $\tilde\bold F$ respectively . Say that second system
{\bf inherits trajectories} of the first system if any trajectory of
the second system as a line (up to a regular reparametrization) coincides
with some trajectory of the first system.
\enddefinition
\definition{Definition 13.3} Two Newtonian dynamical systems on the
Riemannian manifold $M$ are called {\bf trajectory equivalent} if they
inherit trajectories of each other.
\enddefinition
    Note that definitions~13.2 and 13.3 are somewhat different from
corresponding definitions used in papers \cite{54--60}. Our definitions
are more specialized and adopted to the case of Newtonian dynamical
systems with common configuration space.\par
    Metrizable dynamical systems are trivial regarding to their use
in the construction of normal shift. Normal shift along trajectories
of such systems, in essential, is reduced to geodesic normal shift.
Therefore in papers \cite{5} and \cite{9} we considered the problem
of describing all metrizable Newtonian dynamical systems admitting
the normal shift, and in paper \cite{15} we constructed examples of
non-metrizable ones. Main result of papers \cite{5} and \cite{9} is
formulated in the following theorem.
\proclaim{Theorem 13.1} Newtonian dynamical system admitting the normal
shift on the Riemannian manifold $M$ with metric $\bold g$ is metrizable
by means of conformally equivalent metric $\tilde\bold g=e^{-2f}\,\bold
g$ if and only if its force field is given by formula
$$
F_k=-|\bold v|^2\,\nabla_kf+2\sum^n_{s=1}
\nabla_sf\,v^s\,v_k+\frac{H(v\,e^{-f})\,e^f}{|\bold v|}\,v_k,
$$
where $H=H(v)$ is some arbitrary function of one variable.
\endproclaim
    Theorem~13.1 solved the problem of describing dynamical systems
admitting the normal shift and being metrizable by means of conformally
equivalent metric. However, the requirement of conformal equivalence
of metrics $\tilde\bold g$ and $\bold g$ in this theorem is
{\bf a priori}. One can exclude this requirement. Then geodesic flow of
metric $\tilde\bold g$ will correspond to the dynamical system with
less special force field 
$$
F^k=\sum^n_{i=1}\sum^n_{j=1}(\Gamma^k_{ij}-\tilde\Gamma^k_{ij})\,v^i
\,v^j\hskip -3em
\tag13.3
$$
in the metric $\bold g$. One can step further, i\.~e\. one can avoid
metric $\tilde\bold g$ at all assuming $\tilde\Gamma^k_{ij}$ in
\thetag{13.3} to be components of some symmetric affine connection
in $M$.
\definition{Definition 13.4} Newtonian dynamical system on Riemannian
manifold $M$ is called {\bf metrizable by geodesic flow} of affine
connection $\tilde\Gamma$ in $M$ if it inherits trajectories of this
geodesic flow.
\enddefinition
    This is the very treatment of the concept of metrizability that was
considered in paper \cite{11}, which is not published unfortunately. In
paper \cite{11} the following theorem was proved.
\proclaim{Theorem 13.2} Newtonian dynamical system admitting the normal
shift on the Riemannian manifold $M$ of the dimension $n\geqslant 3$ is
metrizable by geodesic flow of affine connection $\tilde \Gamma$ if and
only if this connection is (at least locally) a metric connection
for some metric $\tilde\bold g=e^{-2f}\,\bold g$, which is conformally
equivalent to basic metric $\bold g$ of the manifold $M$.
\endproclaim
    In other words, if there exists some Newtonian dynamical system
which inherits trajectories of geodesic flow of affine connection
$\tilde\Gamma$, and which is admitting the normal shift, then connection
$\tilde\Gamma$ is necessarily a metric connection for some metric
$\tilde\bold g=e^{-2f}\,\bold g$. And conversely, if $\tilde\Gamma$ is
defined by metric $\tilde\bold g=e^{-2f}\,\bold g$, then one can find
some dynamical system admitting the normal shift and inheriting
trajectories of geodesic flow of $\tilde\Gamma$. Though the converse
proposition of the theorem is obvious, since geodesic flow of metric
$\tilde\bold g=e^{-2f}\,\bold g$ admits the normal shift in metric
$\bold g$.\par
    Theorem~13.2 shows that a priori assumption on conformal equivalence
of metrics $\tilde\bold g$ and $\bold g$ used in papers \cite{5} and
\cite{9} at the first approach to the problem of metrizability doesn't
cause the loss of generality. Proof of the theorem~13.2 is based on the
following fact, which was proved in \cite{5}.
\proclaim{Theorem 13.3} Suppose that force field $\bold F$ of the first
Newtonian dynamical system on Riemannian manifold $M$ is a homogeneous
function of degree $2$ with respect to components of velocity vector
$\bold v$ in the fibers of tangent bundle $TM$. Then second Newtonian
dynamical system inherits trajectories of the first system if and only
if its force field $\tilde\bold F$ is given by the following formula:
$$
\tilde\bold F(p,\bold v)=\bold F(p,\bold v)+\frac{H(p,\bold v)}
{|\bold v|}\cdot\bold v.
$$
\endproclaim
    Components of force field \thetag{13.3} are quadratic functions
with respect to the components of velocity vector. Hence the for force
field of dynamical system inheriting trajectories of geodesic flow of
affine connection $\tilde\Gamma$ we have the formula
$$
F^k=\sum^n_{i=1}\sum^n_{j=1}(\Gamma^k_{ij}-\tilde\Gamma^k_{ij})\,v^i
\,v^j+H\,N^k.\hskip -3em
\tag13.4
$$
Further proof of theorem~13.2 in unpublished paper \cite{11} consisted
in substituting \thetag{13.4} into the normality equations \thetag{3.1}
and \thetag{3.2}. Here we shall give more simple proof of this theorem
based on comparison of formulas \thetag{13.4} and \thetag{8.2}. Let's
denote $M^k_{ij}=\tilde\Gamma^k_{ij}-\Gamma^k_{ij}$. The quantities
$M^k_{ij}$ are the components of some (not extended) tensor field
$\bold M$ on the manifold $M$. It is called the field of {\bf deformation}
or the field of {\bf variation} for connection $\Gamma$. Now formula
\thetag{13.4} is written as follows:
$$
F^k=-\sum^n_{i=1}\sum^n_{j=1}M^k_{ij}\,v^i\,v^j+H\,N^k.\hskip -3em
\tag13.5
$$
\demo{Proof of the theorem~13.2} Let's contract both sides of formula
\thetag{13.5} with components of orthogonal projector $\bold P$ from
\thetag{4.3}. This immediately excludes the term containing scalar
function $H=H(x^1,\ldots,x^n,v^1,\ldots,v^n)$:
$$
\sum^n_{k=1}P^q_k\,F^k=-\sum^n_{k=1}\sum^n_{i=1}\sum^n_{j=1}
P^q_k\,M^k_{ij}\,v^i\,v^j.\hskip -3em
\tag13.6
$$
Similar contracting in formula \thetag{8.2} cancels the entry of
function $h(W)$:
$$
\sum^n_{k=1}P^q_k\,F^k=-|\bold v|^2\sum^n_{k=1}P^q_k\,U^k(x^1,\ldots,x^n,W).
\hskip -3em
\tag13.7
$$
Here in \thetag{13.7} by $U^k=U^k(x^1,\ldots,x^n,|\bold v|)$ we denote
the following quantities:
$$
U^k=\sum^n_{i=1}\frac{g^{ki}\,V_i(x^1,\ldots,x^n,W(x^1,\ldots,x^n,
|\bold v|))}{|\bold v|}.
\hskip -3em
\tag13.8
$$
From \thetag{3.10} we see that $U^k$ are components of extended
vector field $\bold U$, which is fiberwise spherically symmetric
(see definition~5.1 above). Let's compare the relationships
\thetag{13.7} and \thetag{13.6}. Their left hand sides coincide.
Hence we can equate right hand sides of these two relationships:
$$
\sum^n_{k=1}\sum^n_{i=1}\sum^n_{j=1}P^q_k\,M^k_{ij}\,v^i\,v^j=
|\bold v|^2\sum^n_{k=1}P^q_k\,U^k.\hskip -3em
\tag13.9
$$
This relationship is remarkable, since the dependence on $\bold v$
in it is almost explicit. Indeed, the quantities $M^k_{ij}$ do not
depend on $\bold v$, while quantities $U^k$ depend only on modulus
of velocity vector. Let's fix the point $p\in M$. This means that
we fix local coordinates $x^1,\,\ldots,\,x^n$. Then rewrite the
relationship \thetag{13.9} in vectorial form:
$$
\bold P(\bold M(\bold v,\bold v)-|\bold v|^2\cdot\bold U)=0.
\hskip -3em
\tag13.10
$$
Here $\bold M(\bold v,\bold v)$ is vector valued quadratic form determined
by tensor $\bold M$ when twice contracting it with vector $\bold v$. Note
that according to \thetag{13.10} the operator of projection $\bold P$, when
applied to the expression $\bold M(\bold v,\bold v)-|\bold v|^2\cdot\bold
U$, yields zero. Therefore the equality \thetag{13.10} can be rewritten as
follows:
$$
\bold M(\bold v,\bold v)-|\bold v|^2\cdot\bold U(|\bold v|)
=\lambda(\bold v)\cdot\bold v.
\hskip -3em
\tag13.11
$$\par
    Vector $\bold U=\bold U(|\bold v|)$ in \thetag{13.11} depend only on
modulus of velocity vector, while scalar $\lambda=\lambda(\bold v)$ can
contain full scale dependence on $\bold v$. Let's study this dependence.
Consider vectors $\bold U(|\bold v|)$ and $\bold U(|\alpha\cdot\bold v|)$,
where $\alpha$ is a number. Remember that in theorem~13.2 we deal with
multidimensional case $n\geqslant 3$. In the space of the dimension
$n\geqslant 3$ vector $\bold v$ ran be turned so that it doesn't belong
to the linear span of vectors $\bold U(|\bold v|)$ and $\bold U(|\alpha
\cdot\bold v|)$, while its modulus $|\bold v|$ being preserved unchanged.
Let's substitute $\alpha\cdot\bold v$ for vector $\bold v$ into the
equality \thetag{13.11}:
$$
\alpha^2\cdot\bold M(\bold v,\bold v)-\alpha^2\,|\bold v|^2\cdot
\bold U(|\alpha\cdot\bold v|)=\alpha\,\lambda(\alpha\cdot\bold v)\cdot
\bold v.\hskip -3em
\tag13.12
$$
Then multiply both sides of \thetag{13.11} by $\alpha^2$ and subtract
the obtained equality from \thetag{13.12}. As a result we get the
following relationship:
$$
\alpha\,|\bold v|^2\cdot(\bold U(|\alpha\cdot\bold v|)
-\bold U(|\bold v|))+(\lambda(\alpha\cdot\bold v)-\alpha\,
\lambda(\bold v))\cdot\bold v=0.\hskip -3em
\tag13.13
$$
For the case when vector $\bold v$ doesn't belong to linear span
of vectors $\bold U(|\bold v|)$ and $\bold U(|\alpha\cdot\bold v|)$
from the equality \thetag{13.13} we derive
$$
\align
&\bold U(|\alpha\cdot\bold v|)=\bold U(|\bold v|),\hskip -3em
\tag13.14\\
\vspace{1ex}
&\lambda(\alpha\cdot\bold v)=\alpha\,\lambda(\bold v).\hskip -3em
\tag13.15
\endalign
$$
Though the equality \thetag{13.14} holds for the case when $\bold v$
belongs to linear span of $\bold U(|\bold v|)$ and $\bold U(|\alpha
\cdot\bold v|)$ as well, since $\bold U$ depend on $|\bold v|$, but
not on the direction of $\bold v$. Substituting \thetag{3.14} back to
\thetag{3.13}, we prove \thetag{3.15} for all $\bold v\neq 0$. For
$\bold v=0$ the value of $\lambda(\bold v)$ is not determined by
formula \thetag{13.11}. Therefore we can extend the function $\lambda(
\bold v)$ by taking $\lambda(0)=0$. This cancels the restriction
$\bold v\neq 0$ in applying formula \thetag{13.15}.\par
     Due to \thetag{13.14} vector $\bold U$ do not depend on $\bold v$
at all. Therefore left hand side of \thetag{13.11} is quadratic function
in $\bold v$. The equality \thetag{13.11} can be rewritten as
$$
\bold K(\bold v)=\bold K(\bold v,\bold v)=\lambda(\bold v)\cdot
\bold v.\hskip -3em
\tag13.16
$$
Any quadratic in $\bold v$ function satisfies the following identity,
\pagebreak which can be checked by direct calculations: $\bold K(\bold
v_1+\bold v_2)+\bold K(\bold v_1-\bold v_2)=2\cdot\bold K(\bold v_1)
+2\cdot\bold K(\bold v_2)$. Substituting \thetag{13.16} into this
identity, we get the equality for $\lambda(\bold v)$:
$$
\lambda(\bold v_1+\bold v_2)\cdot(\bold v_1+\bold v_2)+
\lambda(\bold v_1-\bold v_2)\cdot(\bold v_1-\bold v_2)=
2\,\lambda(\bold v_1)\cdot\bold v_1+2\,\lambda(\bold v_2)
\cdot\bold v_2.
$$
If vectors $\bold v_1$ and $\bold v_2$ are linearly independent, then
the above vectorial equality leads to the pair of scalar equalities:
$$
\align
&\lambda(\bold v_1+\bold v_2)+\lambda(\bold v_1-\bold v_2)=2\,
\lambda(\bold v_1),\\
&\lambda(\bold v_1+\bold v_2)-\lambda(\bold v_1-\bold v_2)=2\,
\lambda(\bold v_2).
\endalign
$$
Let's add them and divide the result by $2$. Then we get the relationship
$$
\lambda(\bold v_1+\bold v_2)=\lambda(\bold v_1)+\lambda(\bold v_2).
\hskip -3em
\tag13.17
$$
If vectors $\bold v_1$ and $\bold v_2$ are linearly dependent, then
\thetag{13.17} follows from \thetag{13.15}. The relationships
\thetag{13.15} and \thetag{13.17} mean that $\lambda(\bold v)$ is a
linear function in $\bold v$.\par
    Thus $\lambda(\bold v)=(\bold v\,|\,\boldsymbol\Lambda)$, where
$\boldsymbol\Lambda$ is some (not extended) vector field on $M$.
When applied to quadratic form $\bold M(\bold v,\bold v)$ in
\thetag{13.11}, this yields
$$
\bold M(\bold v,\bold v)=|\bold v|^2\cdot\bold U+(\bold v\,|\,
\boldsymbol\Lambda)\cdot\bold v.\hskip -3em
\tag13.18
$$
Above we have proved that $\bold U$ doesn't depend on $\bold v$.
Therefore $\bold U$ is also some (not extended) vector field on
$M$. Now let's return to \thetag{13.8} and rewrite this equality in
terms of covariant components of the vector $\bold U$:
$$
U_i(x^1,\ldots,x^n)\,|\bold v|=V_i(x^1,\ldots,x^n,W(x^1,\ldots,x^n,
|\bold v|)).
$$
The quantity $|\bold v|$ plays the role of independent variable in this
equality. Let's substitute $V(x^1,\ldots,x^n,w)$ for $|\bold v|$ and take
into account the relationships \thetag{7.12}:
$$
U_i(x^1,\ldots,x^n)\,V(x^1,\ldots,x^n,w)=V_i(x^1,\ldots,x^n,w).
$$
Now remember the relationships \thetag{7.13}. They show that
$U_i(x^1,\ldots,x^n)$ is a logarithmic derivative of the function
$V(x^1,\ldots,x^n,w)$ with respect to the variable $x^i$. This
logarithmic derivative doesn't depend on $w$:
$$
U_i(x^1,\ldots,x^n)=\frac{\partial\ln(V(x^1,\ldots,x^n,w))}{\partial x^i}.
\hskip -3em
\tag13.19
$$
From \thetag{13.19} it follows that there exist (at least locally)
two functions, a function $f=f(x^1,\ldots,x^n)$ and a function
$\rho=\rho(w)$ such that 
$$
\pagebreak
\aligned
&U_i=\nabla_if=\frac{\partial f(x^1,\ldots,x^n)}{\partial x^i},\\
\vspace{1ex}
&V=\exp(f(x^1,\ldots,x^n))\,\rho(w).
\endaligned\hskip -3em
\tag13.20
$$
At the expense of gauge transformation \thetag{8.3} function $V=V(x^1,
\ldots,x^n,w)$ in \thetag{13.20} can be brought to the following form:
$$
V=\exp(f(x^1,\ldots,x^n))\,w.\hskip -3em
\tag13.21
$$
Let's substitute \thetag{13.21} into the formula \thetag{8.2} for the
components of force field $\bold F$:
$$
F_k=h(|\bold v|\,e^{-f})\,e^f\,N_k+\sum^n_{i=1}|\bold v|^2\,U_i\,
\bigl(2\,N^i\,N_k-\delta^i_k\bigr).\hskip -3em
\tag13.22
$$
Now let's compare formula \thetag{13.22} with formula \thetag{13.5} and
take into account the above relationship \thetag{13.18} for tensor field
$\bold M$. This yields
$$
\frac{H-h(|\bold v|\,e^{-f})\,e^f}{|\bold v|^2}=
(2\,\bold U+\boldsymbol\Lambda\,|\,\bold N).\hskip -3em
\tag13.23
$$
Here $\bold N$ is unitary vector directed along the vector of velocity.
Note that left hand side of \thetag{13.23} depend only on modulus of
velocity vector $\bold v$, while right hand side depend on the direction
of this vector. Therefore from \thetag{13.23} we get
$$
\xalignat 2
&\quad H=h(|\bold v|\,e^{-f})\,e^f,&&\boldsymbol\Lambda=-2\,\bold U.
\hskip -3em
\tag13.24
\endxalignat
$$
The relationships \thetag{13.20} and \thetag{13.24} completely determine
the components of tensor field $\bold M$. For components of connection
$\tilde\Gamma$ we have
$$
\tilde\Gamma^k_{ij}=\Gamma^k_{ij}-\nabla_if\,\delta^k_j
-\nabla_jf\,\delta^k_i+\sum^n_{q=1}g^{kq}\,\nabla_qf\,g_{ij}.
\hskip -3em
\tag13.25
$$
Substituting \thetag{13.25} into \thetag{13.3}, we come to the force
field \thetag{13.2}. Force field \thetag{13.2} corresponds to geodesic
flow of metric connection for metric $\tilde\bold g=e^{-2f}\,\bold g$.
Thus, theorem~13.2 is completely proved.\qed\enddemo
\head
14. Acknowledgments.
\endhead
     Author is grateful to A.~T.~Fomenko and V.~S.~Vladimirov for the
opportunity to report results of thesis \cite{17} in seminars at Moscow
State University and at Steklov Mathematical Institute. Author is
grateful to participants of these seminars for the attention and fruitful
discussion. Author is especially grateful to A.~V.~Bolsinov and
A.~T.~Fomenko for several questions, which appears very stimulating. The
answer to one of these questions is obtained in this paper (see above).
Author also is grateful to A.~V.~Bolsinov for help in arranging contacts
with some geometers from Moscow and Saint-Petersburg.\par
     Work is supported by grant from Russian Fund for Basic Research
(project No\nolinebreak\.~00\nolinebreak-01-00068, coordinator 
Ya\.~T.~Sultanaev), and by
grant from Academy of Sciences of the Republic Bashkortostan
(coordinator N.~M.~Asadullin).\par

\enddocument
\end